\documentclass[11pt, a4paper,leqno]{amsart}
\usepackage{amsmath,amsthm,amscd,amssymb,amsfonts, amsbsy}
\usepackage{latexsym}
\usepackage{exscale}
\usepackage{mathrsfs}
\usepackage{enumitem}
\usepackage[utf8]{inputenc}
\usepackage{datetime}

\usepackage[colorlinks,citecolor=red,pagebackref,hypertexnames=false, breaklinks]{hyperref}

\calclayout
\allowdisplaybreaks

\newdateformat{monthdayyeardate}{%
	\monthname[\THEMONTH]~\THEDAY, \THEYEAR}
\newcommand{\X}{\mathbf{X}}

\newcommand{\normH}[1]{{\vert\kern-0.25ex\vert\kern-0.25ex\vert #1 
		\vert\kern-0.25ex\vert\kern-0.25ex\vert}}

\newcommand{\normHB}[1]{{\Big\vert\kern-0.25ex\Big\vert\kern-0.25ex\Big\vert #1 
		\Big\vert\kern-0.25ex\Big\vert\kern-0.25ex\Big\vert}}

\def\Xint#1{\mathchoice
    {\XXint\displaystyle\textstyle{#1}}%
    {\XXint\textstyle\scriptstyle{#1}}%
    {\XXint\scriptstyle\scriptscriptstyle{#1}}%
    {\XXint\scriptscriptstyle\scriptscriptstyle{#1}}%
    \!\int}
    \def\XXint#1#2#3{{\setbox0=\hbox{$#1{#2#3}{\int}$}
    \vcenter{\hbox{$#2#3$}}\kern-.5\wd0}}
    \def\fint{\Xint-}
    
  \def\Xint#1{\mathchoice
    {\XXint\displaystyle\textstyle{#1}}%
    {\XXint\textstyle\scriptstyle{#1}}%
    {\XXint\scriptstyle\scriptscriptstyle{#1}}%
    {\XXint\scriptscriptstyle\scriptscriptstyle{#1}}%
    \!\int}
    \def\XXint#1#2#3{{\setbox0=\hbox{$#1{#2#3}{\int}$}
    \vcenter{\hbox{$#2#3$}}\kern-.5\wd0}}

\def \W{ \mathbf{W}}

\def \R{ \mathbb{R}}
\def \Z{ \mathbf{Z}}
\def \Y{ \mathbf{Y}}

\def\Pcal{\mathcal{P} }

\theoremstyle{plain}
\newtheorem{theorem}[equation]{Theorem}
\newtheorem{lemma}[equation]{Lemma}

\newtheorem{proposition}[equation]{Proposition}

\newtheorem{claim}[equation]{Claim}

\theoremstyle{definition}
\newtheorem{definition}[equation]{Definition}

\theoremstyle{remark}
\newtheorem{remark}[equation]{Remark}

\newcommand{\Rd}{\mathbb{R}^d_+}

\numberwithin{equation}{section}

\begin{document}
\allowdisplaybreaks
\author{Gael Diebou Yomgne}
\address{Gael Diebou Yomgne \\
 Mathematisches Institut, Rheinische Friedrich-Wilhelms-Universit\"{a}t Bonn\\
Endenicher Allee 60\\
53115, Bonn, Germany }
\email{gaeldieb@math.uni-bonn.de}

\author{Herbert Koch}
\address{Herbert Koch\\
 Mathematisches Institut, Rheinische Friedrich-Wilhelms-Universit\"{a}t Bonn\\
Endenicher Allee 60\\
53115, Bonn, Germany} \email{koch@math.uni-bonn.de}

\title[Dirichlet problem for weakly harmonic maps]{Dirichlet problem for weakly harmonic maps with rough data}
\date{\today}

\thanks{ G.D.Y. was supported by the DAAD through the program \textquotedblleft Graduate School Scholarship Programme, 2018\textquotedblright (Number 57395813) and by the Hausdorff Center for Mathematics at Bonn. H.K. was supported by Deutsche Forschungsgemeinschaft through CRC 1060 and the Hausdorff Center for Mathematics at Bonn}

\keywords{Weakly harmonic maps, small BMO data, strictly stable, Carleson measure}

\begin{abstract}
Weakly harmonic maps from a domain $\Omega$ (the upper half-space $\Rd$ or a bounded $C^{1,\alpha}$ domain, $\alpha\in (0,1]$) into a smooth closed manifold are studied. Prescribing small Dirichlet data in either of the classes $L^{\infty}(\partial\Omega)$ or  $BMO(\partial\Omega)$, we establish solvability of the resulting boundary value problems by means of a nonvariational method. As a by-product, solutions are shown to be locally smooth, $C^{\infty}_{loc}$. Moreover, we show that boundary data can be chosen large in the underlying topologies if $\Omega$ is smooth and bounded by perturbing strictly stable smooth harmonic maps.
\end{abstract}

\maketitle

\tableofcontents

\bigskip
\section{Introduction and statements of main results}\label{S1}
%%%%%%%%%%%%%%%%%%%%%%%%%%%%%%%%%%%%%%%%%%%%%%%%%
In this paper, we are interested in the solvability of the Dirichlet problem for the weakly harmonic maps equation subject to irregular data at the boundary.  Assume for simplicity sake that $\Omega=\mathbb{R}^d_{+}:=\{x=(x',x_d), \hspace{0.1cm}x'\in \R^{d-1},\hspace{0.1cm}x_d>0\}$, ($d>2$) (note that one can also consider $\Omega\subset \R^{d}$ bounded $C^{1,\alpha}$ domain) and let $N$ be a smooth closed  Riemannian manifold. As a result of Nash's embedding theorem we can assume without any restriction that $N$ isometrically embeds into $\mathbb{R}^{m}$ for some positive integer $m$. Denote by $\dot{W}^{1,2}(\mathbb{R}^{d}_{+},N)$ the space of functions $u:\mathbb{R}^d_+\rightarrow \mathbb{R}^m$ whose first order (distributional) derivatives belong to the Lebesgue space $L^2(\mathbb{R}^d_+,\mathbb{R}^{dm})$ and satisfy the constraint $u(x)\in \textit{N}$ a.e. $x\in \mathbb{R}^{d}_{+}$. With $x_1,...,x_{d}$ representing a coordinate system on $\mathbb{R}^{d}_{+}$, one can associate to any Sobolev map $u\in \dot{W}^{1,2}(\mathbb{R}^{d}_{+},N)$ an energy density defined as  \begin{equation*}
e(u)=\displaystyle\sum_{i=1}^{d} \partial_{ x_i}u\cdot\partial_{ x_i}u =\left|\nabla u \right|^2. 
\end{equation*}The Dirichlet energy functional of $u$ is then given by 
\begin{equation*}
E(u)=\displaystyle\dfrac{1}{2}\int_{\mathbb{R}^d_{+}}\textit{e}(u)dx.    
\end{equation*} Consider a tubular neighborhood $U$ of the manifold $N$ in $\R^m$ on which the nearest point projection map $\Pcal_{N}:U\rightarrow N$ is well-defined and smooth (cf. Section \ref{re} below). For any test function $\phi\in C^{\infty}_0(\mathbb{R}^{d}_{+},\mathbb{R}^{m})$ and for $s>0$ small enough, critical points of the functional $E$ are maps $u$ in $\dot{W}^{1,2}(\mathbb{R}^{d}_{+},N)$ such that the first variation of the energy satisfies the identity 
\begin{equation*}
\dfrac{\partial }{\partial s}\bigg|_{s=0}E(\mathcal{P}_{N}(u+s\phi))=0    
\end{equation*} where $\mathcal{P}_{N}(u+s\phi)$ belongs to $\dot{W}^{1,2}(\mathbb{R}^{d}_{+}, N)$. The Euler-Lagrange system associated to this variational problem reads 
\begin{align}\label{E}\Delta u+\Gamma(u)(\nabla u, \nabla u)=0 
\end{align} in the sense of distributions where $\Delta$ is the Laplace operator for the $d$-dimensional Euclidean space and $\Gamma(q):T_qN\times T_qN\rightarrow (T_qN)^{\perp}$ is the second fundamental form of the embedding  $N\hookrightarrow\mathbb{R}^m$ with $$\Gamma(u)(\nabla u, \nabla u)=\displaystyle\sum_{i=1}^{d}\Gamma(u)\left( \dfrac{\partial u}{\partial x_i},\dfrac{\partial u}{\partial x_i}\right).$$ Solutions of Syst. \eqref{E} are called weakly harmonic maps. Supplementing this equation with the boundary condition 
\begin{equation}\label{Dir-data}
u=f\hspace{0.2cm} \mbox{on} \hspace{0.2cm}\partial\R^{d}_+,    
\end{equation}
we aim in this manuscript at addressing the well-posedness issue for the boundary value problem \eqref{E}-\eqref{Dir-data} when $f$ assumes minimal translation and scale invariant regularity assumption in the sense made precise at \eqref{scaling}. The question of existence of  harmonic maps plays an important role in differential geometry, Teichm\"{u}ller theory \cite{J} and in surface matching problem (in computer vision). In hydrodynamics theory, weakly harmonic maps into the $2$-sphere are fundamental objects in the modelling of flows of nematic liquid crystals. 

The Dirichlet problem for weakly harmonic maps in various geometrical settings have been studied in many works. When the source manifold is a compact connected Riemannian manifold of class $C^3$, Hildebrandt and his collaborators \cite{HKW} showed the existence of small solutions in the energy space provided the image of the boundary is contained in a small geodesic ball of $N$. Moreover, they further  proved that if the radius of the ball is strictly bounded above by $\frac{\pi}{2\kappa^{1/2}}$ where $\kappa\geq 0$ is an upper-bound for the sectional curvature of $N$, then solutions are $C^2$-regular in the interior. Uniqueness of these small solutions was independently obtained by J\"{a}ger and Kaul in \cite{JK1}. When $\Omega$ is the unit Euclidean ball in $\R^3$, Struwe \cite{Str} established solvability in the Sobolev class $H^{1,2}(\Omega,N)$ for data $f$ having small energy. In particular, only solutions obeying the restriction 
\begin{equation}\label{smallness-bd-energy}
\sup_{x\in \Omega, r>0}\bigg(r^{-1}\int_{B_r(x)\cap \Omega}|\nabla u|^2dy\bigg)<\varepsilon
\end{equation}
for $\varepsilon>0$ sufficiently small are unique. This result was generalized to arbitrary dimension ($d\geq 3$) in \cite{M}. The study of rotationally symmetric weakly harmonic maps with finite energy and their stability is the main subject of the article \cite{JK}. Regarding regularity, observe that the nonlinearity in Eq. \eqref{E} belongs to $L^1$ whenever $u$ has finite energy. Thus a boostrapping argument will not improve the initial regularity of the solution.
We quote Helein's unconditional regularity results for two-dimensional sources and general targets \cite{H,H1} which use the special structure of the equation in 2D. See also \cite{J,S} for the higher regularity of continuous weakly harmonic maps.

We observe that most of the aforementioned solvability results employ direct methods of calculus of variations under higher regularity condition on boundary data. In order to lower the requirement in smoothness, one needs new techniques. A starting point is to identify a befitting notion of solutions allowing for low regularity data. Formally, the Dirichlet problem \eqref{E}-\eqref{Dir-data} can be reformulated using  Green identities so that the resulting equation reads 
\begin{equation}\label{int-E}
u(x)=\mathcal{H}f(x)+\mathcal{N}(\Gamma(u)(\nabla u,\nabla u))(x),\quad u(x)\in N \hspace{0.1cm} a.e. \hspace{0.1cm}x\in \Rd    
\end{equation}
where $\mathcal{H}f$ is the Poisson extension of $f$ and $\mathcal{N}$ the Newtonian potential, see Section \ref{section2} below for more details. Observe that $v=\mathcal{H}f$ makes sense as an absolutely convergent integral under the weaker condition 
\begin{equation}\label{cond-on-f}
I:=\displaystyle\int_{\mathbb{R}^{d-1}}|f(x')|(1+|x'|^{d})^{-1}dx'<\infty
\end{equation}    
and $v$ solves the Laplace equation in $\Rd$. If $\mathcal{M}$ denotes the centered Hardy-Littlewood maximal function, then it can be verified that $I\leq C\mathcal{M}f(z)$ for some constant $C>0$ depending on a fixed point $z\in \R^{d-1}$. Thus, whether or not \eqref{cond-on-f} holds can be verified using the mapping properties of $\mathcal{M}$. It is worth pointing out that the latter fully characterizes the solvability of the Dirichlet problem for linear elliptic systems of second-order with constant complex coefficients in half-space \cite{MMMM1}. Moreover, if $u$ formally solves Eq. \eqref{E} then, the rescaled map 
\begin{equation}\label{scaling}
 u_{\lambda}(x):=u(\lambda x), \hspace{0.2cm}x\in \Rd,\hspace{0.1cm}\lambda>0    
\end{equation} is another solution since the second fundamental form $\Gamma(u)(\cdot,\cdot)$ has a quadratic growth in the gradient of $u$. Hence, we seek for classes of functions defined on $\mathbb{R}^{d-1}$, enjoying both \eqref{cond-on-f} and the scaling law \eqref{scaling} and such that there exists a suitable notion of trace associated to harmonic functions in half-space. Natural candidates include the John-Nirenberg's space $BMO(\R^{d-1})$ and the class of measurable essentially bounded functions on $\R^{d-1}$. Indeed, it was established by Fefferman \cite{F} that $f\in BMO(\R^{d-1})$ if its Poisson extension $v=\mathcal{H}f(x)$ satisfies 
\begin{equation}
\sup_{x'\in \R^{d-1},x_d>0}\left(x_d^{-(d-1)} \int_{B_{x_d}(x')}\int_0^{x_d}s|\nabla v(y,s)|^2dyds\right)^{1/2}<\infty.
\end{equation}
This condition was later shown to characterize all harmonic functions whose trace belong to $BMO(\R^{d-1})$, see \cite{FJN,FS}.
This motivates the consideration of the following functional setting. 
\begin{definition}\label{def1} Let $d>2$ and $m>1$.
Call $\X$ the space of functions $u:\mathbb{R}^{d}_{+}\rightarrow \R^m$ such that 
\begin{equation}\label{norm-X}
\|u\|_{\X}=\|u\|_{L^{\infty}(\mathbb{R}^d_+)}+|u|_{\X} <\infty    
\end{equation}
\end{definition}
where the semi-norm $|\cdot|_{\X}$ reads 
\begin{equation*}
|u|_{\X}=\displaystyle\sup_{x_d>0}x_d\|\nabla u\|_{L^{\infty}(\R^{d-1})}+\displaystyle\sup_{(x',x_d)\in \Rd}\left(x_d^{1-d}\int_{\displaystyle\textit{B}_{x_d}(x')}\int_{0}^{x_d}y_d|\nabla u|^2dy_ddy' \right)^{1/2} .    
\end{equation*}
When endowed with the norm \eqref{norm-X}, $\X$ is a Banach space and one easily verifies that it is scaling invariant with respect to \eqref{scaling}. Such functional frameworks turn out to be suitable for the analysis of certain critical boundary value problems and their use certainly goes beyong the context of harmonic maps. Our motivation comes from Koch \& Tataru's work \cite{KT} on the well-posedness for the Navier-Stokes equations where similar spaces with parabolic scaling and ideas leading to their consideration were first introduced. This approach was subsequently employed in many other works  \cite{GL,KL1,KL2,W}, just to mention a few. While these articles exclusively deal with parabolic problems, the present work also aims at showing how elliptic boundary value problems subject to low regularity data can be analyzed via similar methods. Observe that the first term in Eq. \eqref{int-E} is harmonic and has a well-defined trace at the boundary while the second term is continuous up to the boundary when $u\in \X$. In what follows, the boundary value problem should be understood in this sense. 

We are ready to state our main results.
\begin{theorem}\label{theo1} Assume that  $f\in L^{\infty}(\mathbb{R}^{d-1},\textit{N})$. There exists a positive number $\varepsilon:=\varepsilon(d,N)$ such that if $\|f\|_{L^{\infty}(\mathbb{R}^{d-1})}\leq \varepsilon$, then the Dirichlet problem \eqref{E}-\eqref{Dir-data} is uniquely solvable in a small closed  ball of $\X$. 
\end{theorem}
The \textit{BMO}-Dirichlet problem for (\ref{E}) is equally well-posed in the following sense.   
\begin{theorem}\label{theo2}
There exists $\varepsilon_0>0$ with the property that for any map $f\in BMO(\mathbb{R}^{d-1},N)$ with $\|f\|_{\textit{BMO}(\mathbb{R}^{d-1})}\leq \varepsilon_0$, the BMO-Dirichlet problem \eqref{E}-\eqref{Dir-data}  admits a solution $u\in \X$. Furthermore, this solution is unique in a small closed ball of $\X$, $B^{\X}_{c\varepsilon}(v)=\{u\in \X:\|u-v\|_{\X}\leq c\varepsilon_0\}$ for some constant $c>0$ depending on $N$. Here $v$ is the harmonic extension of $f$.
\end{theorem}
This existence result seems sharp in the sense that  $BMO(\R^{d-1})$ is the largest translation and scaling invariant (with respect to \eqref{scaling}) space so that the first iteration of the fixed point map is well-defined.   
\begin{remark}
Although our solvability results have been stated for $\Rd$, analogous conclusions for the geometrical setting given by bounded smooth domains remain valid in natural analogues of $\X$ (see Theorem \ref{theo42} for more details). Moreover, 
\begin{itemize}
\item Our method also infers the solvability of the problem $-\Delta u=|\nabla u|^2+F$ in $\Omega$ (with $F\in \Y$ with  small norm if $\Omega=\Rd$ or $F\in \Z$ (c.f. Section \ref{section 4}) with small norm if $\Omega$ is bounded) subject to small Dirichlet data in  $L^{\infty}(\partial\Omega)$ or $BMO(\partial\Omega)$.
\item Unlike the results  quoted earlier, our boundary data are allowed to have unbounded energies. However, it is not clear whether the smallness assumptions on the size of the boundary value can be relaxed. For bounded energy data, it is known that uniqueness fails in absence of appropriate smallness condition (like \eqref{smallness-bd-energy} for weakly harmonic maps from the unit ball in $\R^3$). A similar observation was made in \cite{JK1} for smooth harmonic maps. 
\end{itemize} 
\end{remark}
Despite these evidences about the necessity of having a size restriction condition on the boundary data, one may still ask the question whether or not $f$ in Theorem \ref{theo1} and \ref{theo2} can be prescribed ''large'' in the $L^{\infty}$-norm or BMO semi-norm, respectively. We come back to this particular question in the last part of the manuscript. 
\begin{remark}\label{remark1}
Solutions constructed in the above theorems are locally smooth. Indeed, Theorem \ref{theo1} tells us that $\nabla u$ is bounded in $\mathbb{R}^{d-1}\times(\eta,\infty)$ for any $\eta>0$ since $u\in \X$. Thus $\Delta u\in \textit{L}^{\infty}_{loc}$ and  $u\in W^{2,p}_{loc}$ for any $p<\infty$ by the standard $L^p$-theory for elliptic equations. Once again, by using \eqref{E}, we arrive at $\Delta u \in W^{1,p}_{loc}$ which infers $u \in W^{3,p}_{loc}$. In a repetitive way, one obtains that $u$ belongs to the Sobolev space  $W^{k,p}_{loc}$ for all $k=1,2,3,...$. Applying  Sobolev embedding theorem ultimately yields $u$ in $C_{loc}^{\infty}$. This regularity result, however, is another consequence of smallness -- there are everywhere discontinuous weakly harmonic maps \cite{Ri}. 
\end{remark}

\section{Preliminaries and Auxiliary Results}
\label{section2}
\subsection{The homogeneous theory}
For a point $x$ of $\mathbb{R}^d_+$, we write $x=(x',x_d)$, $x'\in \R^{d-1}$, $x_d\in (0,\infty)$.  Let $f$ be a locally integrable function on $\R^{d-1}$. For a subset $E\subset \R^{d-1}$, denote by $|E|$ its Lebesgue measure and let $f_E=\fint_{E}fdx'=|E|^{-1}\int_{E}fdx'$ the integral mean of $f$. We say that $f$ belongs to $BMO(\R^{d-1})$ if  
 \begin{equation*}
\|f\|_{BMO(\R^{d-1})}=\sup_{B}|B|^{-1}\int_{B}|f(x')-f_{B}|dx'
 \end{equation*} where the supremum is taken over all balls in $\R^{d-1}$. In what follows, we do not always distinguish between $BMO(\R^{d-1},\R^{m})$ and $BMO(\R^{d-1})$ and it should be clear from the context. 
Let  $T(B)= B_r(x')\times (0,r)$ be the Carleson region above the boundary ball $B_r(x')\subset \R^{d-1}$. A measure $\mu$ in $\Rd$ is termed Carleson if  $\mathcal{C}(\mu)=\displaystyle\sup_{B}|B|^{-1}\mu(T(B))$ is finite. The role of these measures in connection to linear elliptic boundary value problems was first observed in \cite{F,FS}: $BMO(\R^{d-1})$ is the trace space of harmonic functions $v$ in $\Rd$ for which $x_d|\nabla v|^2dx'dx_d$ is a Carleson measure. Moreover,  we have the equivalence 
\begin{equation}\label{carleson}
\sup_{x'\in \R^{d-1},x_d>0}\left(x_d^{-(d-1)} \int_{B_{x_d}(x')}\int_0^{x_d}s|\nabla v(y,s)|^2dyds\right)^{1/2}\approx \|f\|_{BMO(\mathbb{R}^{d-1})}.
\end{equation} 

Now consider the Laplace equation $\Delta v=0$ in $\mathbb{R}^d_+$ with $v\big|_{\partial\mathbb{R}^d_+}=f$. The Poisson extension 
\begin{equation}
v(x)=\mathcal{H}f(x):=[P_{x_d}\ast f](x')    
\end{equation} is the unique solution which decays at infinity. Recall that $P_{x_d}$ is explicitly given by 
\begin{align*}
P_{x_d}(x)=x_dP(x/x_d)=c_d\dfrac{x_d}{(|x|^2+x_d^2)^{d/2}}    
\end{align*}
where $c_d$ is a normalizing constant such that $\displaystyle\int_{\R^{d-1}}P_{x_d}(x)dx=1$.
We collect in the following lemma the boundedness properties of $\mathcal{H}$.

\begin{lemma}\label{lem21}
Let $f=(f_1,...,f_m)$ defined on $\mathbb{R}^{d-1}$. Then, $\mathcal{H}f\in \X$ for all $f\in L^{\infty}(\mathbb{R}^{d-1})$ and
\begin{align}\label{l1}\|\mathcal{H}f\|_{\X} \leq C\|f\|_{L^{\infty}(\mathbb{R}^{d-1})}.
\end{align} 
 Moreover, if $f\in BMO(\mathbb{R}^{d-1})$, then there exists a positive constant $C>0$ independent on $f$ such that 
\begin{align}
\label{BMOest}|\mathcal{H}f|_{\X}\leq C\|f\|_{BMO(\mathbb{R}^{d-1})}.
\end{align}  
\end{lemma}
\begin{proof} Note that the estimate \eqref{l1} is invariant by scaling and translation; hence it suffices to consider the case $x'=0$ and $x_d=1$, that is, 
\begin{equation*}
|v(0,1)|+|\nabla v(0,1)|+\|y_d^{1/2}\nabla v\|_{L^2(B_1(0)\times (0,1))}\leq C\|f\|_{L^{\infty}(\mathbb{R}^{d-1})}; \quad v=\mathcal{H}f.
\end{equation*}
Observe that $v$ is harmonic in $\Rd$ so that the bound 
\begin{equation*}
|v(0,1)|+|\nabla v(0,1)|\leq C\|f\|_{L^{\infty}(\mathbb{R}^{d-1})}    
\end{equation*}
follows from local elliptic estimates. To establish the third bound, one proceeds as follows.  Let $B_2(0)$ be a ball in $\mathbb{R}^{d-1}$ and denote by $\chi$ its characteristic function. Decompose $f$ into a local and global part, $f=\chi f+(1-\chi)f=f_1+f_2$. Taking into consideration the  harmonic extension of each part, set $v=v_1+v_2$ and write 
\begin{align*}
\|y_d^{1/2}\nabla v\|_{L^2(B_1(0)\times(0,1))}&\leq \|y_d^{1/2}\nabla v_1\|_{L^2(B_1(0)\times(0,1))}+\|y_d^{1/2}\nabla v_2\|_{L^2(B_1(0)\times(0,1))}\\
&=I_1+I_2\end{align*}
We prove that both $I_1$ and $I_2$ satisfy the desired estimate. For the first integral, we use integration by parts assuming that $f$ is continuous together with the property that $\Delta v=0$ in $\mathbb{R}^{d}_+$. Indeed, 
\begin{align}\nonumber
\|y_d^{1/2}\nabla v\|^2_{L^2(\mathbb{R}^d_+)}\nonumber&=\int_{\mathbb{R}^d_+}y_d|\nabla v|^2dy_ddy'\\
\nonumber&=\sum\limits_{i=1}^{d}\int_{\mathbb{R}^d_+}y_d\partial_iv\cdot \partial_ivdy_ddy'\\
\nonumber&=-\sum\limits_{i=1}^{d}\int_{\mathbb{R}^d_+}(\delta_{id}\partial_iv+y_d\partial^2_iv)\cdot vdy_ddy'\\
\nonumber&=-\int_{\mathbb{R}^d_+}\partial_dv\cdot vdy_ddy'\\
&=\dfrac{1}{2}\int_{\mathbb{R}^{d-1}}f^2(y')dy'.
\end{align}
Here, we can assume without any restriction that $f$ is compactly supported in $B_2(0)$ so that \begin{equation*}
I_1\leq C\|f\|_{L^{\infty}(\mathbb{R}^{d-1})}.    
\end{equation*} On the other hand, using  the kernel decay property 
\begin{align}\label{p4}
|\partial^{\alpha}P(x)|\leq C_n|x|^{1-|\alpha|-d} \hspace{0.2cm}\text{for all}\hspace{0.2cm}x=(x',x_d)\in \overline{\mathbb{R}^d_+}\setminus\{0\}, 
\end{align} 
one obtains
\begin{equation*}
I^2_2=\|y_d^{1/2}\nabla v_2\|^2_{L^2(B_1(0)\times(0,1))}\nonumber=\int_{\mathbb{R}^d_+}y_d|\nabla v_2|^2\chi_{B(0)\times(0,1)}dy_ddy'.    
\end{equation*}	But 
\begin{align*}|\nabla v_2(y',y_d)|&\leq \int|\nabla P^{y_d}(y'-z')||h_2(z')|dz'\\
&\leq \int_{\mathbb{R}^{d-1}\setminus B_2(0)}|\nabla P^{y_d}(y'-z')||h(z')|dz'\\
&\leq C\int_{\mathbb{R}^{d-1}\setminus B_2(0)}(y_d+|y'-z'|)^{-d}|h(z')|dz'.
\end{align*}
Observe that for $y'\in B_1(0)$ and $0\leq y_d\leq1$, it holds that $(y_d+|y'-z'|)^{-d}\leq C(1+|z'|^d)^{-1}$ whenever $z'\in \mathbb{R}^{d-1}\setminus B_2(0)$. Thus 
\begin{align*}
|\nabla v_2(y',y_d)|&\leq C \|f\|_{L^{\infty}(\mathbb{R}^{d-1})}\int_{\mathbb{R}^{d-1}\setminus B_2(0)}(1+|z'|^d)^{-1}dz'\\
&\leq C\|f\|_{L^{\infty}(\mathbb{R}^{d-1})}.
\end{align*}
Upon squaring the previous inequality, multiplying by $y_d$ and integrating over the cylinder $B_1(0)\times (0,1)$ one obtains the bound 
\begin{equation*}
\int_{B_1(0)}\int_{0}^{1}y_d|\nabla v_2(y',y_d)|^2dy_ddy\leq C\|f\|^2_{L^{\infty}(\mathbb{R}^{d-1})}    
\end{equation*}  
which in turn completes the proof of the first statement \eqref{l1}.\\
Next, we prove the estimate \begin{equation}\label{weighted-sup-BMO}
\displaystyle\sup_{x_d>0}x_d\|\nabla v(\cdot,x_d)\|_{L^{\infty}(\mathbb{R}^{d-1})}\leq C\|f\|_{BMO(\mathbb{R}^{d-1})}.
\end{equation}
Since $v$ is harmonic in $\Rd$, then so is the gradient $\nabla v$ and by the mean value theorem, we find that  
\begin{align}\label{mean-value}
|\partial_iv(x)|\leq Cr^{-d}\int_{B_{r}(x)}|\partial_i v|dy,\quad i=\{1,\cdots,d\} 
\end{align} for any ball $B_{r}(x)$ with $\overline{B_r(x)}\subset \Rd$, $r>0$, $x\in \Rd$.  Let $B_{t/4}(x)\subset \mathbb{R}_+^d$ be the ball with center at $x$ and radius $t/4=x_d/3$. It follows that $B_{t/4}(x)\subset Q_t(x')\times [t/2,t]$ where $Q_t(x')$ is the cube in $\R^{d-1}$ with center at $x'$ and side-length $t>0$. We may appeal to \eqref{mean-value} and write 
\begin{align*}
|\partial_iv(x)|&\leq C_d\left( \fint_{B_{t/4}(x)}|\partial_iv|^2dy\right)^{1/2}\\
&\leq C_dt^{-d/2}\left( \int_{Q_t(x')}\int_{t/2}^{t}|\nabla v|^2dy_ddy'\right)^{1/2}\\
&\leq C_dt^{-(d+1)/2}\left( \int_{Q_t(x')}\int_{t/2}^{t}y_d|\nabla v|^2dy_ddy'\right)^{1/2}\\
&\leq C_dt^{-1}\left( t^{1-d}\int_{Q_t(x')}\int_0^{t}y_d|\nabla v(y,y_d)|^2dy_ddy'\right)^{1/2}, 
\end{align*} which in turn implies \eqref{weighted-sup-BMO} and shows \eqref{BMOest} in view of \eqref{carleson}.
\end{proof}
\begin{remark}\label{rmk2}Estimate \eqref{weighted-sup-BMO} can alternatively be derived from the integral representation of $v$ and the cancellation property 
\begin{align}\label{cancelpro}
\int_{\mathbb{R}^{d-1}}\nabla^{\ell}P_{x_d}(x'-y')dy'=0 \hspace{0.2cm}\text{for all}\hspace{0.2cm} (x',x_d)\in \mathbb{R}^d_+, \hspace{0.2cm}\text{for all}\hspace{0.2cm} \ell\in \mathbb{N}.
\end{align}
\end{remark}
\subsection{The Poisson equation}
This section is devoted to the study of the Poisson equation with source term in the space  $\Y$ collecting functions $F$ defined on $\mathbb{R}^d_+$ such that the quantity $\|F\|_{\Y}$ is finite,
\begin{equation*}
\|F\|_{\huge\mathbf{Y}}:=\displaystyle\sup_{x_d>0}x^2_d\| F(\cdot,x_d)\|_{L^{\infty}(\mathbb{R}^{d-1})}+\sup_{(x',x_d)\in \mathbb{R}^d_+}x_d^{1-d}\int_{\displaystyle\textit{B}_{x_d}(x')}\int_{0}^{x_d}y_d| F|dy_ddy'.
\end{equation*}
It is clear that $\Y$ equipped with the above norm is complete, hence Banach. Let $F$ be a measurable function in $\Rd$ such that 
\begin{equation}\label{abs-conv-cond}
\int_{\Rd}\dfrac{y_d|F(y)|}{(1+|y|)^d}dy<\infty.
\end{equation} A solution $u$ to the Poisson equation $-\Delta u=F$ in $\mathbb{R}^d_+$, $u=0$ on $\mathbb{R}^{d-1}$ can  explicitly be given by the Newtonian potential \begin{equation*}u(x)=\mathcal{N}F(x):=\displaystyle\int_{\mathbb{R}^d_+}G(x,y)F(y)dy    
\end{equation*} with $G(\cdot,\cdot)$ the  Green kernel for the Laplacian in the upper half-space $\mathbb{R}^d_+$ which satisfies the following upper-bound estimates   

\begin{enumerate}\item For every $x,y\in \mathbb{R}^d_+$, $x\neq y$ 
\begin{align*}
G(x,y)\leq C\min\left\lbrace \dfrac{\min (x_d,y_d)}{|x-y|^{d-1}},\dfrac{x_dy_d}{|x-y|^{d}},\dfrac{1}{|x-y|^{d-2}}\right\rbrace. 
\end{align*}
\item \label{Green-gradient-bds}For every $x,y\in \mathbb{R}^d_+$, $x\neq y$, \begin{align*}|\nabla G(x,y)|\leq \min\big\{|x-y|^{1-d},y_d|x-y|^{-d}\big\}.
\end{align*}	
\item For each $k\in \mathbb{N}$, \begin{align*}|\nabla^k G(x,y)|\leq |x-y|^{2-k-d}\hspace{0.2cm} \text{for all} \hspace{0.2cm}x,y\in \mathbb{R}^d_+, x\neq y.
\end{align*}
\end{enumerate}
It should be observed that functions in $\Y$ satisfy \eqref{abs-conv-cond}. 
Our next lemma deals with the mapping properties of the Green potential. 
\begin{lemma}\label{lem22}
The Newtonian potential $\mathcal{N}$ maps  $\mathbf{Y}$ boundedly into $\mathbf{X}$ i.e. for any $F\in \mathbf{Y}$, $\mathcal{N}F\in \mathbf{X}$ and in addition, there holds the estimate 
\begin{equation}\label{Newtonian-potential-bds}
\|\mathcal{N}F\|_{\mathbf{X}}\leq C\|F\|_{\mathbf{Y}}
\end{equation}
where the constant $C$ only depends on the dimension d.
\end{lemma}
\begin{proof}Once again, due to the scaling and translation invariance nature of \eqref{Newtonian-potential-bds}, its validity simplifies to that of the following localized bound 
\begin{equation*}
|\mathcal{N}F(0,1)|+|\nabla \mathcal{N}F(0,1)|+\big\|y^{1/2}_d|\nabla\mathcal{N}F|\big\|_{L^2(B_1(0)\times(0,\infty))}\leq C\|F\|_{\mathbf{Y}}\end{equation*}
whose proof is divided into two steps.\\
	
\textbf{Step 1.}  \textit{The inequality} $|\mathcal{N}F(0,1)|\leq C\|F\|_{\mathbf{Y}}$. By definition of $\mathcal{N}$, we have\\
\begin{align*}
|\mathcal{N}F(0,1)|&\leq \int_{\mathbb{R}^{d-1}}\int_{0}^{\infty}G(-y',1-y_d)|F(y)|dy'_ddy\\
&=I_1+I_2+I_3+I_4\end{align*}
where 
\begin{align*}
I_1=\int_{B_1(0)}\int_{0}^{1/2}G(-y',1-y_d)|F(y)|dy_ddy',\hspace{0.2cm} I_2=\int_{B_1(0)}\int_{1/2}^2G(-y',1-y_d)|F(y)|dy_ddy',\\ I_3=\int_{B^c_1(0)}\int_{0}^2G(-y',1-y_d)|F(y)|dy_ddy',\hspace{0.2cm} I_4=\int_{\mathbb{R}^{d-1}}\int_{2}^{\infty}G(-y',1-y_d)|F(y)|dy_ddy'. 	
\end{align*}	
In what follows, we repeatedly make use of the above upper-bound estimates on the Green function $G(\cdot,\cdot)$.	
\begin{align*}
I_1&=\int_{B_1(0)}\int_{0}^{1/2} G(-y',1-y_d)|F(y)|dy_ddy'\\
	&\leq C \int_{B_1(0)}\int_{0}^{1/2}\dfrac{y_d}{(|y'|^2+(1-y_d)^2)^{^{\frac{d-1}{2}}}}|F(y)|dy_ddy'\\
	&\leq C2^{d-1}\int_{B_1(0)}\int_{0}^{1/2}y_d|F(y)|dy_ddy'\\
&\leq C\|F\|_{\Y}. 
\end{align*}
Moving on, we have 
\begin{align*}
I_2&=\int_{B_1(0)}\int_{1/2}^{2} G(-y',1-y_d)|F(y)|dy_ddy'\\
&\leq C \int_{B_1(0)}\int_{1/2}^{2}\dfrac{|F(y)|dy_ddy'}{(|y'|^2+(1-y_d)^2)^{^{\frac{d-2}{2}}}}\\
&\leq C \sup_{y_d>0}y^2_d\|F(\cdot,y_d)\|_{L^{\infty}(\mathbb{R}^{d-1})}\int_{B_1(0)}\int_{1/2}^{2}|y'|^{2-d}dy_{d}dy'\\
&\leq C\|F\|_{\mathbf{Y}}\int_{|y'|\leq 1}\int_{1/2}^{2}|y'|^{2-d}dy_ddy'\\ 
 &\leq C\|F\|_{\mathbf{Y}}.
\end{align*}
To estimate $I_3$, cover $B^c_1(0)=\mathbb{R}^{d-1}\setminus B_1(0)$ with the family of cubes $\left\lbrace Q_1(z')\right\rbrace_{z'\in \mathbb{Z}^{d-1}} $ centered at $z'$, $|z'|>1$ with side length 1. It follows that 
\begin{align*}I_3&\leq \int_{\mathbb{R}^{d-1}\setminus B_1(0)}\int_{0}^{2} G(-y',1-y_d)|F(y)|dy_ddy'\\
&\leq C\sum_{\substack{z'\in \mathbb{Z}^{d-1}\\|z'|>1}}\int_{Q_1(z')\cap B^c_1(0)}\int_{0}^{2}\dfrac{y_d|F(y)|dy_ddy'}{(|y'|^2+(1-y_d)^2)^{^{d/2}}}\\
	&\leq C\sum_{\substack{z'\in \mathbb{Z}^{d-1}\\|z'|>1}}|z'|^{-d}\int_{Q_1(z')}\int_{0}^{2}y_d|F(y)|dy_ddy'\\
	&\leq C\|F\|_{\mathbf{Y}}.
	\end{align*}
Finally, noticing that for all $y_d\geq 2$, $y_d-1\geq y_d/2$ we obtain that 
\begin{align*}
I_4&=\int_{\mathbb{R}^{d-1}}\int_{2}^{\infty} G(-y',1-y_d)|F(y)|dy_ddy'\\
&\leq C \int_{\mathbb{R}^{d-1}}\int_{2}^{\infty}\dfrac{y_d|F(y)|}{(|y'|^2+(y_d-1)^2)^{d/2}}dy_ddy'\\
&\leq C\|y_d^2F\|_{L^{\infty}(\mathbb{R}^d_+)}\int_{\mathbb{R}^{d-1}}\int_{2}^{\infty}\dfrac{dy_ddy'}{y_d[|y'|^2+y_d^2]^{d/2}}\\
&\leq C\|F\|_{\mathbf{Y}}\left( \int_{\mathbb{R}^{d-1}}(1+|z'|^2)^{d/2}dz'\right) \left( \int_{2}^{\infty}y_d^{-2}dy_d\right) \\
&\leq C\|F\|_{\mathbf{Y}}. 
\end{align*}
The strategy employed above via decomposition of the integral domain can also be used to  prove the \textit{local pointwise gradient estimate}. Indeed, the relation  
\begin{equation*}
\nabla\mathcal{N}F(x)=\int_{\Rd}\nabla G(x,y)F(y)dy    
\end{equation*}
holds in the sense of distributions so that by utilizing the pointwise bound \eqref{Green-gradient-bds} and splitting the above integral exactly as before in the same regions to get, say $J_i$, $i=1,2,3,4$, we obtain the desired estimate. We note however that to estimate $J_2$, one rather invokes the property that $|\nabla G(x,\cdot)|\in L^{\frac{d}{d-1},\infty}(\Rd)$ uniformly for any $x\in \Rd$. Here, $L^{p,\infty}(\mathbb{R}^d)$ denotes the Lorentz space defined as the set of measurable functions $f$ such that $\displaystyle \sup_{E}|E|^{1/p-1}\int_{E}|f(y)|dy$ is finite where the supremum is taken over all open subsets of $\R^{d}$. This concludes step 1.\\

\hspace{-0.43cm}\textbf{Step 2.} \textit{The energy estimate.} \\One may proceed here as in the proof of Lemma \ref{lem21} using the Green's kernel bounds together with the usual cut-off procedure on $F$ but there is an alternative shorter argument which allows us to derive this energy-type bound. In fact, there are two different estimates leading to the desired $\textit{L}^2$-bound, namely
\begin{equation*}
\|\mathcal{N}F\|_{\textit{L}^{\infty}(\mathbb{R}^d_+)}\leq C\|F\|_{\mathbf{Y}},
\end{equation*}
whose validity has already been justified in step 1 and the second bound 
\begin{equation*}
\big\|x_d^{1/2}|\nabla \mathcal{N}F|\big\|^2_{L^{2}(\mathbb{R}^d_+)}\leq \|\mathcal{N}F\|_{L^{\infty}(\mathbb{R}^d_+)}\|x_dF\|_{L^1(\mathbb{R}^d_+)}
\end{equation*}
which may be deduced from a priori estimates. We prove the latter estimate by assuming that $F\in \mathbf{Y}$ is smooth and has compact support in $\Rd$. Thus $\mathcal{N}F$ is smooth and it is not difficult to justify the formal calculations below. Now multiply the equation $-\Delta \mathcal{N}F=F$ in $\Rd$ by $x_d\mathcal{N}F$ and integrate by parts over $\mathbb{R}^d_+$ to obtain
\begin{align*}-\int_{\mathbb{R}^d_+}(x_d\mathcal{N}F)\cdot\Delta\mathcal{N}Fdx&=\int_{\mathbb{R}^d_+}x_dF\cdot\mathcal{N}Fdx.\end{align*}
The left hand side of this identity further simplifies to  
\begin{align*}-\int_{\mathbb{R}^d_+}(x_d\mathcal{N}F)\cdot\Delta\mathcal{N}Fdx&=-\sum_{i=1}^{m}\sum_{j=1}^{d}\int_{\mathbb{R}^d_+}\partial_j(x_d\mathcal{N}F_i)\partial_j\mathcal{N}F_idx\\
&=-\sum_{i=1}^{m}\sum_{j=1}^{d}\int_{\mathbb{R}^d_+}\bigg(x_d|\partial_j\mathcal{N}F_i|^2+\partial_jx_d\mathcal{N}F_i\partial_j\mathcal{N}F_i\bigg)dx\\
&=\int_{\mathbb{R}^d_+}x_d|\nabla \mathcal{N}F|^2dx+\sum_{i=1}^{m}\int_{\mathbb{R}^d_+}\partial_d(\mathcal{N}F_i)\mathcal{N}F_idx\\
&=\int_{\mathbb{R}^d_+}x_d|\nabla \mathcal{N}F|^2dx.\end{align*}
As such, with the aid of H\"{o}lder's inequality, this implies \begin{align*}\int_{\mathbb{R}^d_+}x_d|\nabla \mathcal{N}F|^2dx&=\int_{\mathbb{R}^d_+}x_dF\cdot\mathcal{N}Fdx\\
&\leq\|\mathcal{N}F\|_{L^{\infty}(\mathbb{R}^d_+)}\|x_dF\|_{L^1(\mathbb{R}^d_+)}
\end{align*}    
which completes this particular step and finishes the proof of Lemma \ref{lem22}. 
\end{proof}
Next, one shall prove that the  solution $u$ satisfies the pointwise constraint $u(x)\in N$ a.e. $x\in \mathbb{R}^d_+$. A first step towards achieving this is to estimate the distance between the solution and the target manifold $N$ in the $L^{\infty}$-norm.   
\begin{proposition}\label{prop} Assume that $f$ is a measurable bounded map into $N$, compact smooth manifold. There exists a constant $C>0$ independent of $f$ with
\begin{equation}\label{dist-to-N-Linfty}
\mathrm{dist}(v(x),N)\leq C\|f\|_{L^{\infty}(\mathbb{R}^{d-1})}
\end{equation}
for all $x\in\mathbb{R}^d_+$. For all $\varLambda>0$, there exists $C_1>0$ depending on $\varLambda$ and $N$ such that if $f$ belongs to $BMO(\mathbb{R}^{d-1},N)$, then  
\begin{equation}\label{dist-to-N-BMO}
\mathrm{dist}(v(x),N)\leq C_1\|f\|_{BMO(\mathbb{R}^{d-1})}+\varLambda
\end{equation}
for every $x\in\mathbb{R}^d_+$. Here, $v$ is the Poisson extension of $f$. 
\end{proposition}
\begin{proof}We only give the proof of \eqref{dist-to-N-BMO} for the first statement \eqref{dist-to-N-Linfty} directly follows from the fact that the distance function is evaluated with respect to the $\sup$ norm. Pick a real number $\ell>0$, fix $x=(x',x_d)$ in $\mathbb{R}^d_+$ and put $\overline{f_{x}}=\displaystyle\fint_{B_{\ell}(0)}f(x'-x_dy')dy'$. Owing to the triangle inequality, one has the bound $\text{dist}(v(x),N)\leq \text{dist}(v(x',x_d),\overline{f_{x}})+\text{dist}(\overline{f_{x}},N).$ Since $f(y')\in N$ for all $y'\in \mathbb{R}^{d-1}$, we find that $$\text{dist}(\overline{f_{x}},N)\leq |\overline{f_{x}}-f(x'-x_dz')| \hspace{0.2cm}\text{for any} \hspace{0.2cm}z'\in B_{\ell}(0),$$ from which we easily deduce the bound 
\begin{equation}\label{first-BMO-bd}
\text{dist}(\overline{f_{x}},N)\leq \|f\|_{BMO(\mathbb{R}^{d-1})}.
\end{equation} 
Also note that the Poisson kernel for the Laplace operator obeys $P_{x_d}(x')=x_d^{1-d}P(x'/x_d)$ with $|P(x')|\leq\dfrac{C}{(|x'|^2+1)^{d/2}}.$ This permits us to write 
\begin{equation*}
v(x)=\displaystyle\int_{\mathbb{R}^{d-1}}P(y')f(x'-x_dy')dy',
\end{equation*} from which  it follows that \begin{align*}|v(x)-\overline{f_{x}}|&=\bigg|\int_{\mathbb{R}^{d-1}}P(y')[f(x'-x_dy')-\overline{f_{x}}]dy'\bigg|\\
&=\bigg|\bigg( \int_{B_{\ell}(0)}+\int_{\mathbb{R}^{d-1}\setminus B_{\ell}(0)}\bigg)P(y')[f(x'-x_dy')-\overline{f_{x}}]dy' \bigg|\\
&\leq  C_d\int_{B_{\ell}(0)}\dfrac{|f(x'-x_dy')-\overline{f_{x}}|}{(|y'|^2+1)^{d/2}}dy'+2\|f\|_{L^{\infty}(\mathbb{R}^{d-1})}\int_{\mathbb{R}^{d-1}\setminus B_{\ell}(0)}|P(y')|dy'\\
&\leq C_d\ell^{d-1}\|f\|_{BMO}+C\|f\|_{L^{\infty}(\R^{d-1})}\int^{\infty}_{\ell}\dfrac{r^{d-2}}{(1+r^2)^{d/2}}dr.
\end{align*}
With an appropriate choice of the radius $\ell$ depending on $\varLambda$ and $N$, one can achieve \begin{align*}
\displaystyle C\|f\|_{L^{\infty}(\R^{d-1})}\int^{\infty}_{\ell}\dfrac{r^{d-2}}{(1+r^2)^{d/2}}dr\leq \varLambda.    
\end{align*}
Consequently, one has
\begin{equation*}
|v(x)-\overline{f_{x}}|\leq C'\|f\|_{BMO(\mathbb{R}^{d-1})}+\varLambda;\quad C':=C'(\varLambda,N),    
\end{equation*} which combined with \eqref{first-BMO-bd} gives the desired estimate.
\end{proof}

\subsection{Reformulation of the problem}\label{re}

The Dirichlet problem for \eqref{E} coupled with a boundary condition \eqref{Dir-data} can be recast as \begin{align}\label{r1}u(x)=v(x)+\mathcal{N}(\Gamma(u)(\nabla u, \nabla u))(x); \hspace{0.2cm}u(x)\in N \hspace{0.2cm}\text{a.e}\hspace{0.2cm}x\in \mathbb{R}^d_+
\end{align} where $v$ represents the harmonic extension of $f$ and $\mathcal{N}$ the Newtonian potential. Thus, for boundary data which are small in the $L^{\infty}$-norm (and $BMO$-semi-norm) one can uniquely solve Eq. \eqref{r1} via  Banach fixed point argument. However, there is an incompatibility which emanates from the fact that $u$ is thought of as an $\mathbb{R}^m$-valued map whereas the second fundamental form $\Gamma(\cdot)$ must be defined on $N$. To override this, we  construct an extension of $\Gamma$ to the entire space $\mathbb{R}^m$. In this respect, of significance to us is the nearest point projection map whose Hessian is expressed in terms of the second fundamental form $\Gamma$.
We know (see for instance \cite[Appendix to chapter $2$, Theorem $1$]{simon}) that if $N$ is a compact smooth manifold isometrically embedded in $\mathbb{R}^m$, then $N$ has a $\rho$-neighborhood in $\mathbb{R}^m$ of the form $U_{\rho}=\{z\in \mathbb{R}^m: \text{dist}(z,N)<\rho\}$ such that the projection $\mathcal{P}_{N}$ which maps a point $z\in U_{\rho}$ to the closest point in $N$ is well-defined and smooth. In addition, it satisfies a number of properties which we partially recall below:
\begin{enumerate}[label={($\textbf{a}_{\arabic*})$}]
\item \label{projectiona1}$\mathcal{P}_{N}(z)\in N,\hspace{0.2cm}z-\mathcal{P}_{N}(z)\in (T_{\mathcal{P}_{N}(z)}N)^{\perp}$, \hspace{0.2cm}$|\mathcal{P}_{N}(z)-z|=\textit{dist}(z,N)$ for all $z\in U_{\rho}$
	\item  $|y-z|>\textit{dist}(z,N)$ for any $y\in N\setminus\{\mathcal{P}_{N}(z)\} $, for all $z\in U_{\rho}$
	\item \label{projection3}$\mathcal{P}_{N}(z+y)=z$ for $z\in N$, $y\in (\textit{T}_zN)^{\perp}$, $|y|<\rho$ and $D_{V}\mathcal{P}_{N}|_z=\textit{P}^{\perp}_{\mathcal{P}_{N}(z)}(V)$, $z\in U_{\rho}$, $V\in \mathbb{R}^m$
	\item\label{projectiona4} $\textit{Hess}\hspace{0.1cm}\mathcal{P}_{N}(z)(V_1,V_2)=-\Gamma(z)(V_1,V_2)$ for $z\in N$ and $V_1,V_2\in \textit{T}_z N$  
\end{enumerate}  
where $D_V$ stands for the directional derivative in the direction of $V$, $\textit{P}^{\perp}_{\mathcal{P}_{N}(z)}$ denotes the orthogonal projection of $\mathbb{R}^m$ onto $\textit{T}_{\mathcal{P}_{N}(z)}N$ and $\textit{Hess}\hspace{0.1cm}\mathcal{P}_{N}(z)$ denotes the Hessian of $\mathcal{P}_{N}$ at $z$.\\
We then extend the second fundamental form $\Gamma$ as follows: take a smooth extension of the projection $\mathcal{P}_{N}$, say $\mathcal{P}\in C^{\infty}(\mathbb{R}^m,\mathbb{R}^m)$ so that $\mathcal{P}$ restricted to $U_{\rho}$ coincides with $\mathcal{P}_{N}$ and define the extension $\widetilde{\Gamma}$ of $\Gamma$ by \begin{align*}\widetilde{\Gamma}(z)(V,W)=-\textit{Hess}\hspace{0.1cm}\mathcal{P}(z)(V,W), \hspace{0.1cm}z\in \mathbb{R}^m;\hspace{0.2cm} V,W\in T_z\mathbb{R}^m.\end{align*}
In the special case $N=\mathbb{S}^{m-1}$, the nearest point projection map $\mathcal{P}_{\mathbb{S}^{m-1}}$ can be realized explicitly. Clearly, the set $U_{\frac{1}{4}}=\{u\in\mathbb{R}^m: \frac{3}{4}\leq |u|\leq \frac{5}{4}\}$ is a neighborhood of $\mathbb{S}^{m-1}$ in $\mathbb{R}^m$ and one may consider $\mathcal{P}_{\mathbb{S}^{m-1}}:U_{\frac{1}{4}}\rightarrow \mathbb{S}^{m-1}$, $u\mapsto \dfrac{u}{|u|}$. 
Now, introduce the operator $\mathcal{S}$ defined  by   \begin{align}\label{fixpteq}\mathcal{S}u(x)=v(x)+\mathcal{N}[\widetilde{\Gamma}(u)(\nabla u,\nabla u)](x),\hspace{0.3cm} x\in \mathbb{R}^d_+.
\end{align}
Note that the above formulation contains the information at boundary. 
Our problem then becomes that of finding a map $u=(u_1,...,u_m)$ such that 
\begin{align}
\label{fixed-point}u=\mathcal{S}u \hspace{0.2cm}\text{in}\hspace{0.2cm} \mathbb{R}^d_+; \hspace{0.2cm}u\in N\hspace{0.2cm} \text{a.e.} \hspace{0.2cm}\text{in}\hspace{0.1cm}\mathbb{R}^d_+. 
\end{align} 
In the sequel, we study some basic properties of the operator $\mathcal{S}$, especially  those required for an eventual application of the Banach fixed point theorem. To this end, since $v\in \textbf{X}$ thanks to Lemma \ref{lem21} consider the closed ball $B_{\epsilon}^{\mathbf{X}}(v)\subset \mathbf{X}$ centered at $v$ with radius $\varepsilon$, \begin{align*}B_{\varepsilon}^{\mathbf{X}}(v)=\{u\in \mathbf{X}: \|u-v\|_{\mathbf{X}}\leq \varepsilon\}.\end{align*} 
\begin{lemma}\label{lem24}
Assume that the Dirichlet data $f:\mathbb{R}^{d-1}\rightarrow N$ satisfies $\|f\|_{L^{\infty}(\mathbb{R}^{d-1})}\leq \varepsilon$. For all $u$ in the ball $B^{\mathbf{X}}_{\varepsilon}(v)$, there exists $C>0$ depending only on $d$ with $\|u\|_{\mathbf{X}}\leq C\varepsilon.$
\end{lemma}  
The proof of this result  immediately follows from  Lemma \ref{lem21}. Likewise, we have
\begin{lemma}\label{bmolem1}
Given $f\in BMO(\mathbb{R}^d_+;N)$ with $\|f\|_{BMO(\mathbb{R}^{d-1})^m}\leq \varepsilon$, the following estimates hold, namely 
\begin{align*}|u|_{\mathbf{X}}\leq c\varepsilon,\hspace{0.2cm} \|u\|_{L^{\infty}(\mathbb{R}^d_+)}\leq C \hspace{0.2cm}\text{for all}\hspace{0.2cm}u\in \textit{B}^{\mathbf{X}}_{\varepsilon}(v)
\end{align*} where $c:=c(d)>0$ and $C:=C(\varepsilon,N)>0$.
\end{lemma}
The next result establishes the mapping properties of $\mathcal{S}$ and arises as a direct consequence of the Lemmas \ref{lem24} and \ref{bmolem1}.
\begin{lemma}\label{lem25}
There exists $\varepsilon'>0$ such that the operator $\mathcal{S}$ maps $B^{\mathbf{X}}_{\varepsilon'}(v)$ into itself whenever $f\in L^{\infty}(\mathbb{R}^{d-1})$ satisfies the smallness condition $\|f\|_{L^{\infty}(\mathbb{R}^{d-1})}\leq \varepsilon'$. Furthermore,  if $\|f\|_{BMO(\mathbb{R}^{d-1})}\leq\varepsilon'$ then $\mathcal{S}:B^{\mathbf{X}}_{\varepsilon'}(v)\rightarrow B^{\mathbf{X}}_{\varepsilon'}(v)$ continuously with respect to the semi-norm $|\cdot|_{\X}$.
\end{lemma} 
\begin{proof}Let $\varepsilon_0>0$ and let $u\in B^{\mathbf{X}}_{\varepsilon_0}(v)$, we want to achieve $\|\mathcal{S}u-v\|_{\mathbf{X}}\leq \varepsilon_0$. Taking into account \eqref{fixpteq} and by using the potential estimate from Lemma \ref{lem22}, it follows that \begin{align*}\|\mathcal{S}u-v\|_{\mathbf{X}}&=\|\mathcal{N}\widetilde{\Gamma}(u)(\nabla u,\nabla u)\|_{\mathbf{X}}\leq C\|\widetilde{\Gamma}(u)(\nabla u,\nabla u)\|_{\mathbf{Y}}\\
&=C\bigg(\sup_{x_d>0}x_d^2\|\widetilde{\Gamma}(u)(\nabla u,\nabla u)(\cdot,x_d)\|_{L^{\infty}(\mathbb{R}^{d-1})}+\\
&\hspace{2.15cm}\qquad{}\sup_{(x',x_d)\in\mathbb{R}^d_+}x_d^{1-d}\int_{B_{x_d}(x')}\int_{0}^{x_d}y_d|\widetilde{\Gamma}(u)(\nabla u,\nabla u)(y',y_d)|dy_ddy' \bigg)\\
&\leq C\bigg(\sup_{x_d>0}x_d^2\|\nabla u(\cdot,x_d)\|^2_{L^{\infty}(\mathbb{R}^{d-1})}+\\
&\hspace{3.75cm}\qquad{}\sup_{(x',x_d)\in\mathbb{R}^d_+}x_d^{1-d}\int_{B_{x_d}(x')}\int_{0}^{x_d}y_d|\nabla u(y',y_d)|^2dy_ddy' \bigg)\\
&\leq C\bigg(\|x_d\nabla u(\cdot,x_d)\|_{L^{\infty}(\mathbb{R}^{d}_+)}+\\
&\hspace{2.5cm}\qquad{}\sup_{(x',x_d)\in\mathbb{R}^d_+}\bigg(x_d^{1-d}\int_{B_{x_d}(x')}\int_{0}^{x_d}y_d|\nabla u(y',y_d)|^2dy_ddy'\bigg)^{1/2}\bigg)^2\\
\|\mathcal{S}u-v\|_{\mathbf{X}}&\leq C\|u\|^2_{\mathbf{X}}\leq C\varepsilon_0^2\leq \varepsilon_0
\end{align*}
due to Lemma \ref{lem24} as long as $\varepsilon_0$ is chosen small enough. Mimicking the preceding proof, one obtains the second part of the Lemma (relying in this case on Lemma \ref{bmolem1}) whenever $f$ is sufficiently small in $BMO$.   
\end{proof}

\begin{lemma}\label{lem26}
Let $\varepsilon'>0$ be as in Lemma \ref{lem25}. There exists $\varepsilon_0\in (0,\varepsilon')$  with the property that if $\|f\|_{L^{\infty}(\mathbb{R}^{d-1})}\leq \varepsilon_0$ then the operator $\mathcal{S}:B^{\mathbf{X}}_{\varepsilon_0}(v)\rightarrow B^{\mathbf{X}}_{\varepsilon_0}(v)$ is a contraction map, that is, there exists $\theta\in(0,1)$ with 
\begin{align*}
\|\mathcal{S}u-\mathcal{S}w\|_{\mathbf{X}}\leq\theta\|u-w\|_{\mathbf{X}} \hspace{0.2cm}\text{for all}\hspace{0.2cm}u,w\in\textit{B}^{\mathbf{X}}_{\varepsilon_0}(v).
\end{align*} 
\end{lemma}
\begin{proof}
By linearity of $\mathcal{N}$ and in light of Lemma \ref{lem22}, we have that \begin{align*}\|\mathcal{S}u-\mathcal{S}w\|_{\mathbf{X}}&=\big\|\mathcal{N}[\widetilde{\Gamma}(u)(\nabla u, \nabla u)-\widetilde{\Gamma}(w)(\nabla w, \nabla w)]\big\|_{\mathbf{X}}\\
&\leq C \big\|(\widetilde{\Gamma}(u)(\nabla u, \nabla u)-\widetilde{\Gamma}(u)(\nabla w, \nabla w))+\\
&\hspace{6cm}(\widetilde{\Gamma}(u)(\nabla w, \nabla w)-\widetilde{\Gamma}(w)(\nabla w, \nabla w))\big\|_{\Y}\\
&\leq C(J_1+J_2)
\end{align*} where $$J_1=\big\||\nabla (u- w)| ( |\nabla u|+|\nabla w|)\big\|_{\mathbf{Y}}\hspace{0.2cm}\text{and}\hspace{0.2cm} J_2=\big\||u- w||\nabla w|^2 \big\|_{\mathbf{Y}}.$$ 
We estimate $J_1$ using H\"{o}lder's inequality and Lemma \ref{lem24}.
\begin{align*}J_1&=\sup_{x_d>0}x_d^2\||\nabla (u- w)| ( |\nabla u|+|\nabla w|)(\cdot,t)\|_{L^{\infty}(\mathbb{R}^{d-1})}+\\
&\hspace{2.73cm}\qquad{}\sup_{x',x_d>0}x_d^{1-d}\int_{B_{x_d}(x')}\int_0^{x_d}y_d|\nabla (u- w)| ( |\nabla u|+|\nabla w|)dy_ddy'\\
&\leq\displaystyle\sup_{x_d>0}x_d\|\nabla (u- w) \|_{L^{\infty}(\mathbb{R}^{d-1})}\sup_{x_d>0}x_d\big(\|\nabla u \|_{L^{\infty}(\mathbb{R}^{d-1})}+\|\nabla w \|_{L^{\infty}(\mathbb{R}^{d-1})}\big)+\\
&\qquad{}\hspace{3cm}\sup_{x',x_d>0}\bigg(x_d^{1-d}\int_{B_{x_d}(x')}\int_0^{x_d}y_d(|\nabla u|+|\nabla w|)^2dy_ddy'\bigg)^{1/2}\times\\
\qquad{}&\hspace{4.8cm}\sup_{x',x_d>0}\bigg(x_d^{1-d}\int_{B_{x_d}(x')}\int_0^{x_d}y_d|\nabla (u- w)|^2dy_ddy'\bigg)^{1/2}\\
J_1&\leq C(|u|_{\X}+|w|_{\X})\bigg(\displaystyle\sup_{x_d>0}x_d\|\nabla (u- w) \|_{L^{\infty}(\R^{d-1})}+\\
\qquad{}&\hspace{4.3cm}\sup_{x',x_d>0}\bigg(x_d^{1-d}\int_{B_{x_d}(x')}\int_{0}^{x_d}y_d|\nabla (u- w)|^2 dy_ddy'\bigg)^{1/2}\bigg)\\
J_1&\leq C\varepsilon_0\bigg(\displaystyle\sup_{x_d>0}x_d\|\nabla (u- w) \|_{L^{\infty}(\mathbb{R}^{d-1})}+\\
&\hspace{3.4cm}\qquad{}\sup_{x',x_d>0}\bigg(x_d^{1-d}\int_{B_{x_d}(x')}\int_{0}^{x_d}y_d|\nabla (u- w)|^2 dy_ddy'\bigg)^{1/2}\bigg)\\
J_1&\leq C\varepsilon_0\|u-w\|_{\X}.
\end{align*}
Estimating $J_2$ does require the use of the inequality $a^2+b^2\leq (a+b)^2$, $a,b\geq0$ and Lemma \ref{lem24}. 
\begin{align*}J_2&=\big\||\nabla w|^2|u-w|\big\|_{\mathbf{Y}}\\
&=\sup_{x_d>0}x_d^2\||\nabla w|^2|u-w|\|_{L^{\infty}(\mathbb{R}^{d-1})}+\sup_{x',x_d>0}x_d^{1-d}\int_{B_{x_d}(x')}\int_0^{x_d}y_d|\nabla w|^2| u- w|dy_ddy'\\
&\leq \sup_{x_d>0}x_d^2\|\nabla w\|^2_{L^{\infty}(\mathbb{R}^{d-1})}\|u-w\|_{L^{\infty}(\mathbb{R}^{d}_+)}+\\
&\hspace{4.1cm}\qquad{}\|u-w\|_{L^{\infty}(\mathbb{R}^{d}_+)}\sup_{x',x_d>0}x_d^{1-d}\int_{B_{x_d}(x')}\int_0^{x_d}y_d|\nabla w|^2dy_ddy'\\
&\leq \|u-w\|_{L^{\infty}(\mathbb{R}^{d}_+)}\bigg(\sup_{x_d>0}x_d\|\nabla w\|_{L^{\infty}(\mathbb{R}^{d-1})}+\\
&\hspace{5.1cm}\qquad{}\sup_{x',x_d>0}\bigg(x_d^{1-d}\int_{B_{x_d}(x')}\int_0^{x_d}y_d|\nabla w|^2dy_ddy'\bigg)^{1/2}\bigg)^2\\
&\leq C \|u-w\|_{L^{\infty}(\mathbb{R}^{d}_+)}\|w\|^2_{\X}\\
J_2&\leq C\varepsilon^2_0 \|u-w\|_{\mathbf{X}}.
	\end{align*} 
Summarizing, we find that $$\|\mathcal{S}u-\mathcal{S}w\|_{\mathbf{X}}\leq C\varepsilon_0(1+\varepsilon_0)\|u-w\|_{\mathbf{X}}.$$ One can make $\theta=C\varepsilon_0(1+\varepsilon_0)<1$ if $\varepsilon'$ is chosen sufficiently small. This achieves the proof of Lemma \ref{lem26}.\end{proof}
Following the lines of the above proof, we easily deduce the following.
\begin{lemma}\label{lembmo2}
Let $\varepsilon'>0$ be the number in Lemma \ref{lem25}. There exists $\varepsilon_1\in (0,\varepsilon')$ and $\theta_0\in(0,1)$ such that whenever $f\in BMO(\mathbb{R}^{d-1})$ with $f(x')\in N$ a.e. $x'\in \mathbb{R}^{d-1}$ satisfies $\|f\|_{BMO(\mathbb{R}^{d-1})}\leq \varepsilon_1$, the operator $\mathcal{S}:B^{\X}_{\varepsilon_1}(v)\rightarrow B^{\X}_{\varepsilon_1}(v)$ is a $\theta_0$-contraction map, that is, $$\|\mathcal{S}u-\mathcal{S}w\|_{\X}\leq\theta_0\|u-w\|_{\X} \hspace{0.2cm}\text{for all}\hspace{0.2cm}u,w\in B^{\X}_{\varepsilon_1}(v).$$ 
\end{lemma}

\section{Proofs of the main results}\label{S3}
This section aims at proving Theorems \ref{theo1} and \ref{theo2} by making use of the auxiliary results derived in the previous section.
\subsection{Proof of Theorem \ref{theo1}} A simple application of the contraction principle establishes the existence and uniqueness of solutions. The main task is to show that the solution $u$ satisfies the constraint $u\in N$. 
\begin{proof}In light of Lemmas \ref{lem25} and \ref{lem26} and the Banach fixed-point Theorem, there exists $\varepsilon^0:=\varepsilon^0({N,d})$ such that for  $\|f\|_{L^{\infty}(\mathbb{R}^{d-1})}\leq \varepsilon^0$, Eq. \eqref{fixed-point} admits a unique small solution in $\X$. Now, we need to show that lies in $N$ using  Proposition \ref{prop}. As announced in Section 2, we first show that the distance from the solution $u$ to $N$ can be appropriately controlled so that $u$ lives in a tubular neighborhood of $N$. Let $x\in \mathbb{R}^d_+$, by virtue of Proposition \ref{prop}, we have that
\begin{align*}\text{dist}(u(x),N)&\leq \text{dist}(v(x),N)+\sup_{x\in\mathbb{R}^d_+}|u(x)-v(x)|\\&\leq c\|f\|_{L^{\infty}(\mathbb{R}^{d-1})}+\|u-v\|_{L^{\infty}(\mathbb{R}^{d}_+)}\\
&\leq c\varepsilon_0+\|u-v\|_{\X}\\
&\leq C\varepsilon_0.
\end{align*}  
This implies that $u\in U_{\rho}$ for so long as $c\varepsilon_0<\rho$. As a consequence, one obtains the following identity  
\begin{align}\label{iden}\Delta u= \nabla^2\mathcal{P}_{N}(u)(\nabla u, \nabla u) \hspace{0.2cm}\text{in}\hspace{0.2cm}\mathbb{R}^d_+,\end{align} from \ref{projectiona4} (see section $2.3$). Now define for $z\in U_{\rho}$, the map $\Upsilon_{N}(z)=z-\mathcal{P}_{N}(z)$ and observe that the conclusion immediately follows if $\Upsilon_{N}(u)$ vanishes identically. The existence theory reveals that the gradient of the solution $u$ to the Dirichlet problem \eqref{E}-\eqref{Dir-data} is locally bounded; this qualitative property  as pointed out in Remark \ref{remark1} entails higher regularity of $u$. In fact, $u\in C_{loc}^{\infty}\cap L^{\infty}(\Rd)$, $\Upsilon_{N}(u)$ is bounded and the following holds in the weak sense
\begin{align*}-\Delta \bigg(\dfrac{1}{2}|\Upsilon_{N}(u)|^2\bigg)
&=-\left\langle  \Upsilon_{N}(u),\Delta \Upsilon_N(u) \right\rangle-|\nabla\Upsilon_{N}(u)|^2\\
&=-\left\langle \Upsilon_{N}(u),\nabla^2\Upsilon_{N}(u)(\nabla u, \nabla u)\right\rangle-\left\langle \Upsilon_{N}(u),\nabla\Upsilon_N(u)(\Delta u)\right\rangle-|\nabla\Upsilon_{N}(u)|^2\\
&=\left\langle \Upsilon_{N}(u),\nabla^2\mathcal{P}_{N}(u)(\nabla u, \nabla u)\right\rangle+\left\langle\Upsilon_{N}(u),\nabla\Upsilon_N(u)(\Delta u)\right\rangle-|\nabla\Upsilon_{N}(u)|^2\\
&=-|\nabla\Upsilon_{N}(u)|^2
\end{align*}   
where we have successively used besides the formulas $\nabla \Upsilon_{N}(z)(p)=(\textit{Id}-\nabla \mathcal{P}_{N}(z))(p)$ and $\nabla^2 \Upsilon_{N}(z)(p,q)=-\nabla^2 \mathcal{P}_{N}(z)(p,q)$ for all $p,q\in \mathbb{R}^m$, the identity \eqref{iden} together with the properties \ref{projectiona1} and \ref{projection3}. Indeed, it holds that \[\Upsilon_{N}(u)\in (T_{\mathcal{P}_{N}(u)}N)^{\perp}\hspace{0.3cm} \mbox{and}\hspace{0.3cm} \nabla\mathcal{P}_{N}(u)(\nabla^2\mathcal{P}_{N}(u)(\nabla u, \nabla u))\in T_{\mathcal{P}_{N}(u)}N\] for $u\in U_{\rho}$. On the other hand, since $\mathcal{P}_{N}$ coincides with the identity map of $N$ at the boundary, it follows that $\Upsilon_{N}(u)=0$ on $\mathbb{R}^{d-1}$. Hence, $G(u)=\frac{1}{2}|\Upsilon_{N}(u)|^2$ is a bounded subharmonic function in $\Rd$, one can apply the maximum principle (see e.g. \cite{BCN}) to obtain $u=\mathcal{P}_{N}(u)\in N$.  The proof of Theorem \ref{theo1} is now complete.\end{proof}

\subsection{Proof of Theorem \ref{theo2}}
Here, we argue similarly as before given that the auxiliary results used in the proof of Theorem \ref{theo1} have analogous versions for data sitting in $BMO(\mathbb{R}^{d-1})$.        
\begin{proof}Thanks to Lemmas \ref{lem25} and \ref{lembmo2}, an application of the Banach fixed point Theorem shows that Eq. \eqref{fixed-point} has a unique solution $u$ in $\mathcal{B}^{\X}_{C\varepsilon_1}(0)=\{w\in \X:\hspace{0.1cm} |w|_{\X}\leq C\varepsilon_1\}$  for some constant $C>0$ whenever $f$ satisfies the smallness condition $\|f\|_{BMO(\mathbb{R}^{d-1})}\leq \varepsilon_1$.  In the next lines, we prove that $u \in N$. In effect, it follows from Proposition \ref{prop} that $\text{dist}(u(x),N)\leq \text{dist}(v(x),N)+\varepsilon_1\leq C_1\varepsilon_1+C_2$ for any $x\in \mathbb{R}^d_+$. This shows in particular that $u\in U_{\rho_0}$ provided $C_1\varepsilon_1+C_2< \rho_0$. Therefore, $\Upsilon_{N}(u)\big|_{\R^{d-1}}=0$ and by a similar argument to that performed above we conclude that $u(x)\in N$ a.e. $x\in \mathbb{R}^d_+$.
\end{proof}

\section{Large data situation in bounded domains}\label{section 4}
In this section we attempt to answer the question whether or not the Dirichlet problem for weakly harmonic maps equations subject to ``large'' data in $BMO(\R^{d-1})$ or $L^{\infty}(\R^{d-1})$ is solvable (in the sense described in Section \ref{S1}). In such a scenario, it is clear, based on the theory which has been developed earlier that the norm of the solution in our function space $\X$ may grow, leading to a nonexistence result. This motivates the consideration of stable smooth solutions of equation (\ref{E}) and more specifically, boundary data which are to a certain sense close to the latter -- we shall be more precise regarding this statement in subsequent lines. Existence of stable harmonic maps is not a restricting assumption as exemplified by the class of harmonic maps into targets with nonpositive sectional curvature (see Remark \ref{stabrmk} below). Another class of stable harmonic maps includes local minimizers of the energy $E$ which are smooth under suitable conditions on boundary data. Opting for a perturbation technique the main difficulty comes from the nonlinear geometric constraints in the problem which we bypass by considering an appropriate extension problem and  maximum principle arguments as performed in the proof of Theorems \ref{theo1} and \ref{theo2}. However, we will need the source manifold to be bounded unlike the case treated earlier involving the half-space domain. Consider the weakly harmonic maps equation
 \begin{align}\label{largedataeq}
 -\Delta u=\Gamma(u)(\nabla u,\nabla u)\hspace{0.2cm}\mbox{in}\hspace{0.2cm}\Omega
 \end{align} 
subject to the Dirichlet boundary condition 
\begin{equation}\label{largedata-bc}
u\rvert_{\partial\Omega}=f    
\end{equation}
where $\Omega\subset\mathbb{R}^d$ is a $C^{1,\alpha}$, ($\alpha\in (0,1]$) bounded domain.
Our main Theorems in Section \ref{S1} claims that if $f$ has a small $L^{\infty}$ or $BMO$ norm, then Problem \eqref{largedataeq}-\eqref{largedata-bc} is solvable in $\W$. Before we define the notion of stability, recall that a weak solution to BVP \eqref{largedataeq}-\eqref{largedata-bc} is a map $u\in L^{\infty}\cap W^{1,2}_{f}(\Omega,N)=\{v\in L^{\infty}\cap W^{1,2}(\Omega,\mathbb{R}^m): v\in N\hspace{0.2cm}\mbox{a.e. on}\hspace{0.2cm}\Omega \hspace{0.2cm}\mbox{and}\hspace{0.2cm}v\big|_{\partial\Omega}=f\}$ such that \begin{align}\label{weaksol}\int_{\Omega} \bigg\{ ( \nabla u,\nabla \phi) +\sum_{j=1}^{d}\Gamma(u)(\partial_j u,\partial_j u)\cdot \phi\bigg\}dx=0  
 \end{align} 
 for all $\phi \in L^{\infty}\cap W^{1,2}_{0}(\Omega,\mathbb{R}^m)$
 where $(\cdot,\cdot)$ denotes the standard scalar product in $\mathbb{R}^{dm}$.         
\begin{definition}
Let $u$ be a nontrivial weak solution of Eq. \eqref{largedataeq}. We say that $u$ is strictly stable if 
\begin{align}
\label{stabcond}
Q_{u}(\phi):=\int_{\Omega}\bigg\{ |\nabla \phi|^2+\sum_{j=1}^{d}R^N(\phi,\partial_j u)\partial_j u\cdot \phi\bigg\} dx\geq M\|\phi\|^2_{L^{2}(\Omega)}
\end{align}
for every $\phi\in L^{\infty}\cap W^{1,2}_0(\Omega,\mathbb{R}^m)$, $\phi\in T_uN$ a.e. and for some $M>0$. In case $M=0$, we say that $u$ is stable. 
\end{definition}
In the above definition $R^N(\cdot,\cdot)\cdot$ denotes the curvature tensor of $N$ which at each point $u$ of $N$ is a trilinear map on $T_uN\times T_uN\times T_uN$ to $T_uN$.
The notion of strict stability involving a weighted $L^2$-norm on the R.H.S. of \eqref{stabcond} appeared  in the study of harmonic maps with prescribed set of singularities \cite{H-M}. Note that the integral expression in \eqref{stabcond} up to a change in sign in the second term due to the symmetry feature $(R^N(U,V)W,Z)=-R^N((U,V)Z,W)$ represents the second variation formula of the energy functional associated to \eqref{largedataeq}, see \cite{S}. Thus stability of $u$ and nonnegativeness of the second variation of the energy are formally two equivalent notions.    
\begin{remark}\label{stabrmk}
If the target $N$ has a nonpositive sectional curvature (in the sense of distributions) then any weakly harmonic map $u:\Omega\rightarrow N$ is stable, that is, $u$ satisfies  \eqref{stabcond} with $M=0$. 
\end{remark}

\subsection{Function spaces and linear estimates}
 Recall the definition of the space of bounded mean oscillations defined on the boundary  
 \begin{equation*}
 BMO(\partial\Omega)=\big\{f\in L^1_{loc}(\partial\Omega): \|f\|_{BMO(\partial\Omega)}=\sup_{S\subset \partial\Omega}\sigma(S)^{-1}\int_{S}|f(\xi)-f_{S}|d\sigma(\xi)\big\} 
 \end{equation*} where $S=S_r(\zeta)$, $\zeta\in\partial\Omega$ is the surface ball centered at $\zeta$ with radius $r>0$ and \begin{equation*}f_{S}=\displaystyle\fint_{S}fd\sigma=\dfrac{1}{\sigma(S)}\displaystyle\int_{S}fd\sigma.
 \end{equation*}  
Call  $T(S_r(\zeta))=\Omega\cap B_r(\zeta)$  the Carleson region associated to the surface ball $S_r(\zeta)$ and let $r_0>0$. A measure $\mu$ in $\Omega$ is termed Carleson if there exists a constant $C>0$ depending on $r_0$ such that for all $r\leq r_0$, 
$\mu(T(S_r))\leq C\sigma(S_r)$ with $d\sigma$ being the Lebesgue surface measure.  Fabes \& Neri in  \cite{FN} showed that  harmonic functions on $\Omega$ whose traces belong to $BMO(\partial\Omega)$ can also be characterized by means of Carleson measures. Indeed,  $u$ is harmonic in $\Omega$ and the measure $d\mu(x)=|\nabla u(x)|^2d(x)dx$ ($d(x)$ is the distance from $x$ to $\partial\Omega$) is Carleson if and only if $u$ is the Poisson integral of $f\in BMO(\partial\Omega)$. In addition,  
\begin{align}\label{BMO-bd}
\sup_{S\subset\partial\Omega}\bigg(\sigma(S)^{-1}\int_{T(S)}d(y)|\nabla u(y)|^2dy\bigg)^{1/2} \leq C\|f\|_{BMO(\partial\Omega)}
\end{align}
where the supremum runs over all surface balls $S\subset\partial\Omega$. Next, we introduce some function spaces which will be useful in the sequel. For $\xi\in \partial\Omega$ such that $d(x)=|x-\xi|$, we use the shorthand  $T(S)$ for the Carleson region associated to the surface ball $S_{d(x)}(\xi)$, i.e.  $T(S_{d(x)}(\xi))$.
\begin{definition}
We say that $w:\Omega\rightarrow \mathbb{R}^m$ belongs to $\W$ if the quantity $\|w\|_{\W}$ is finite where 
\begin{align*}
\|w\|_{\W}=\sup_{x\in\Omega}|w(x)|+\sup_{x\in\Omega}d(x)|\nabla w(x)|+\sup_{x\in \Omega}\bigg(d(x)^{1-d}\int_{T(S)}d(y)|\nabla w(y)|^2dy\bigg)^{1/2}.
\end{align*}  
We denote by $\Z$ the space of functions $F:\Omega\rightarrow\mathbb{R}^m$ such that \begin{align*}\|F\|_{\Z}=\sup_{x\in\Omega}d(x)^2|F(x)|+\sup_{x\in \Omega}d(x)^{1-d}\int_{T(S)}d(y)|F(y)|dy<\infty.
\end{align*} 
\end{definition}
 Observe that the function spaces $\W$ and $\Z$ are simply the analogs of $\X$ and $\Y$ in bounded domains respectively, which have been used earlier in Section 1, the fundamental difference being that the distance function to the boundary in this case is a bounded function. Therefore, it is not surprising that some of the results derived in Section 2 persist here. This is the case of the lemma below which provides some relevant information on the solutions to both homogeneous (subject to $BMO$ and $L^{\infty}$ boundary data) and inhomogeneous (with source term in $\Z$) problems for the Laplacian. In what follows we set $[w]_{\W}$ to be   \begin{equation*}
[w]_{\W}:=\sup_{x\in\Omega}d(x)|\nabla w(x)|+\sup_{x\in \Omega}\bigg(d(x)^{1-d}\int_{T(S)}d(y)|\nabla w(y)|^2dy\bigg)^{1/2}.     
 \end{equation*}
\begin{lemma}\label{linearest}
Let $F\in \Z$ and $u$ such that $-\Delta u=F$ in $\Omega$ with $u$$\mid_{\partial\Omega}=f$. The following conclusions hold:
\begin{enumerate}
\item If $f\in BMO(\partial\Omega)$, then $u$ satisfies 
\begin{align*}[u]_{\W}&\leq C(\|f\|_{BMO(\partial\Omega)}+\|F\|_{\Z}).\end{align*}
\item For $f$ in $L^{\infty}(\partial\Omega)$, $u$ is an element of $\W$ and it holds that
\begin{align*}\|u\|_{\W}\leq c(\|f\|_{L^{\infty}(\partial\Omega)}+\|F\|_{\Z}).
\end{align*}
\end{enumerate} 
The generic constants $C$ and $c$ appearing in the above estimates only depend on the dimension and $\Omega$. 
\end{lemma}  
\begin{proof}We distinguish between two steps.\\ 
\hspace{-0.5cm}\textbf{Step 1}. Assume that $F=0$. We prove the corresponding two claims of the lemma.\\ Let $f\in BMO(\partial\Omega,\mathbb{R}^m)$, we would like to establish the bound $[u]_{\W}\leq c\|f\|_{BMO(\partial\Omega)}$. Note that from the Carleson measure characterization of $BMO(\partial\Omega)$, one only needs to verify that the estimate $$\displaystyle\sup_{x\in \Omega}d(x)|\nabla u(x)|\leq c\|f\|_{BMO(\partial\Omega)}$$ is valid. Pick $x_0\in \Omega$, put $R_0=\frac{d(x_0)}{2}$ and assume that $B_{3R_0}(x_0)\subset\subset\Omega$. By harmonicity of $D_iu$ ($i=1,2,\cdots,d$) and standard interior estimates for the Laplace equation we have 
\begin{align*}|\nabla u(x_0)|^2&\leq c\fint_{B_{R_0/4}(x_0)}|\nabla u|^2dy\\
	&\leq cR_0^{1-d}d(x_0)^{-1}\int_{B_{R_0/4}(x_0)}|\nabla u|^2dy.\end{align*}
Let $\xi_0,\xi_y\in \partial\Omega$ such that $d(x_0)=|x_0-\xi_0|$ and $d(y)=|y-\xi_y|$. Since $|y-x_0|\leq d(x_0)/8$, we have $|y-\xi_0|\leq 9R_0/4$ and $7d(x_0)/8\leq d(y)$ so that
\begin{align*}
|\nabla u(x_0)|^2&\leq cR^{1-d}_0d(x_0)^{-2}\int_{T(S_{\frac{9R_0}{4}}(\xi_0))}d(y)|\nabla u|^2dy,
\end{align*}
from which it follows that 
\begin{align*}d(x_0)^{2}|\nabla u(x_0)|^2&\leq c(9R_0/4)^{1-d}\int_{T(S_{9R_0/4}(\xi_0))}d(y)|\nabla u|^2dy.
\end{align*} 
At this point, as $x_0$ was chosen arbitrary we simply pass to the supremum over $\Omega$ on both sides of the above inequality to deduce the desired estimate. For $f$ bounded, the bounds 
\begin{equation}\label{L-infty bd}
\|u\|_{L^{\infty}(\Omega)}\leq c\|f\|_{L^{\infty}(\partial\Omega)}\hspace{0.2cm}\mbox{and}\hspace{0.2cm}\sup_{x\in \Omega}d(x)|\nabla u(x)|\leq C\|f\|_{L^{\infty}(\partial\Omega)}
\end{equation}follows from elliptic interior estimates.  Next, we prove the Carleson measure estimate \begin{align}\label{carlest}\sup_{x\in \Omega}\bigg(d(x)^{1-d}\int_{T(S)}d(y)|\nabla u(y)|^2dy\bigg)^{1/2}\leq C\|f\|_{L^{\infty}(\Omega)}.
\end{align} 
Fix $x\in \Omega$ and let $2S:=S_{2d(x)}(\xi)$ be the surface ball with center at  $\xi\in \partial\Omega$ and radius $2d(x)$. Denoting by $\textbf{1}_{2S}$  the characteristic function of $2S$, we make the decomposition $f=\textbf{1}_{2S}f+(1-\textbf{1}_{2S})f=f_1+f_2$ and write correspondingly $u=u_1+u_2$ for the Poisson extension of $f$ to $\Omega$. We first prove \eqref{carlest} for $u_2$ with the aid of the following pointwise decay bound for the Poisson kernel for the Laplacian in $\Omega$:  \begin{equation*}
|\nabla_x\mathcal{P}(x,\zeta)|\leq \dfrac{cd(x)}{|x-\zeta|^{d+1}},  \quad \zeta\in \partial\Omega    
\end{equation*} which can be found e.g. in \cite{St}. In effect, let $y\in T(S)$, we have  \begin{align*}|\nabla u_2(y)|&\leq \int|\nabla \mathcal{P}(y-\zeta)||f_2(\zeta)|d\sigma(\zeta) \\
&\leq cd(y)\int_{\partial\Omega\setminus(2S)}|y-\zeta|^{-(d+1)}|f(\zeta)|d\sigma(\zeta)\\
&\leq cd(y)\|f\|_{L^{\infty}(\partial\Omega)}\sum_{i=1}^{\infty}\int_{S_i}|y-\zeta|^{-(d+1)}d\sigma(\zeta) 
\end{align*} 
where $S_i=2^{i+1}S\setminus2^{i}S$. Let $\zeta \in S_i$, we have $2^{i}d(x)< |\zeta-\xi|$ which by the triangle inequality implies $|y-\zeta|\geq |\xi-\zeta|-|y-\xi|>2^{i}d(x)-d(x)\geq 2^{i-1}d(x)$. Hence, we have that \begin{align*}|\nabla u_2(y)|&\leq cd(y)d(x)^{-2}\|f\|_{L^{\infty}(\partial\Omega)}\sum_{i=1}^{\infty}2^{(1-i)(d+1)}2^{(i+1)(d-1)}\\
|\nabla u_2(y)|&\leq Cd(x)^{-1}\|f\|_{L^{\infty}(\partial\Omega)}
\end{align*}
since $d(y)\leq d(x)$ whenever $y\in T(S)$. Squaring the above inequality, multiplying both sides by $d(y)$ and integrating over the Carleson region $T(S)$, we obtain \begin{align*}\bigg(d(x)^{1-d}\int_{T(S)}d(y)|\nabla u_2(y)|^2dy\bigg)^{1/2}\leq C\|f\|_{L^{\infty}(\partial\Omega)}.
\end{align*}  This yields the bound we were looking for after passing to the supremum over $\Omega$ on both sides.
To establish the corresponding estimate for $u_1$, we assume without any restriction that $f_1$ is supported in $2S$ and that $\|u_1\|_{L^{\infty}(\Omega)}$ is finite (in view of estimate \eqref{L-infty bd}). Thus, we have  
\begin{align}\label{nontangential-max-bd}
\int_{T(S)}d(y)|\nabla u_1(y)|^2dy\leq C\int_{\Omega}d(y)|\nabla u_1|^2dy\leq C\int_{\partial\Omega}[N^{\ast}u_1(\zeta)]^2d\sigma(\zeta)
\end{align} 
where $N^{\ast}u_1(\zeta)=\displaystyle\sup_{\varGamma(\zeta)}|u_1(x)|$ is the nontangential maximal function of $u_1$, $\varGamma(\zeta)$ is the cone in $\Omega$ with vertex at $\zeta\in \partial\Omega$. We note that the last estimate in \eqref{nontangential-max-bd} is due to Dahlberg \cite{Da}. Hence, from the mapping properties of the nontangential maximal function in Lebesgue space (see \cite{Da1}) and the fact that $\Omega$ is smooth, we find from \eqref{nontangential-max-bd} that 
\begin{align*}
\int_{T(S)}d(y)|\nabla u_1(y)|^2dy\leq C\|f_1\|^2_{L^2(\partial\Omega)}\leq Cd(x)^{d-1}\|f\|^2_{L^{\infty}(\partial\Omega)}.
\end{align*} 
This shows that \eqref{carlest} is valid and completes this part.\\
\textbf{Step 2.}  We prove that any solution $u$ of $-\Delta u=F$ in $\Omega$ which vanishes on $\partial\Omega$ satisfies the bound 
\begin{align}\|u\|_{\W}\leq C\|F\|_{\Z}.
\end{align}
Under the condition that $\Omega$ satisfies the uniform exterior sphere condition, the Green function $G_{\Omega}$ for $\Delta$  satisfies (see \cite[Theorem 3.3]{GW})
\begin{align*}
G_{\Omega}(x,y)\leq C\min\left( \dfrac{1}{|x-y|^{d-2}},\dfrac{d(y)}{|x-y|^{d-1}},
\dfrac{d(x)d(y)}{|x-y|^d}\right)\end{align*} and 
\begin{align*}
|\nabla G_{\Omega}(x,y)|\leq C\min\left( |x-y|^{1-d},\dfrac{d(y)}{|x-y|^{d}}\right).
\end{align*} 
To derive the $L^{\infty}$-estimate, we may write $u$ as the Green potential of $F$ such that for $x\in \Omega$, we have 
\begin{align*}
|u(x)|&\leq \int_{\Omega}|G_{\Omega}(x-y)||F(y)|dy\\
&\leq \left\lbrace  \int_{\{y\in \Omega:\hspace{0.1cm}|x-y|\leq 2^{-1}d(x)\}}+\int_{\{y\in \Omega:\hspace{0.1cm}|x-y|>2^{-1}d(x)\}}\right\rbrace |G_{\Omega}(x-y)||F(y)|dy\\
&=I+II. 
\end{align*}
Since for any $y\in B_{2^{-1}d(x)}(x)$ we have the inequality $d(y)\geq d(x)/2$, we handle $I$ as follows 
\begin{align}\label{firstest}
\nonumber I&\leq c\sup_{x\in \Omega}(d^2(y)|F(y)|)\int_{\Omega\cap B_{2^{-1}d(x)}(x)}d^{-2}(y)|x-y|^{2-d}dy\\
\nonumber &\leq c\|F\|_{\Z}d^{-2}(x)\int_{B_{2^{-1}d(x)}(x)}|x-y|^{2-d}dy\\
I&\leq c\|F\|_{\Z}.
\end{align} In order to estimate the second integral, we cover the set $\big\{y\in \Omega:\hspace{0.1cm}|x-y|>2^{-1}d(x)\big\}$ with the family of annuli $(A_i)_{i}$, $A_i=2^{i}B_{d(x)}(x)\setminus 2^{i-1}B_{d(x)}(x)$ and use the above pointwise estimate on the Green kernel to arrive at 
\begin{align*}
II&\leq \sum_{i=0}^{\infty}\int_{A_i}|F(y)|G_{\Omega}(x-y)dy\\
&\leq C\sum_{i=0}^{\infty}\int_{A_i}\dfrac{d(x)d(y)|F(y)|}{|x-y|^{-d}}dy\\
&\leq Cd^{1-d}(x)\sum_{i=0}^{\infty}2^{-(i-1)d}\int_{2^{i}B_{d(x)}(x)}d(y)|F(y)|dy.
\end{align*}
It is easy to see that $y\in T(2^{i+1}S)$ whenever $y\in 2^{i}B_{d(x)}(x)$ so that
\begin{align*}
II &\leq \sum_{i=0}^{\infty}2^{-(i-1)d}2^{(i+1)(d-1)} (2^{i+1}d(x))^{1-d}\int_{T(2^{i+1}S)}d(y)|F(y)|dy\leq C\|F\|_{\Z}. 
\end{align*}
Combining the latter with \eqref{firstest} yields the desired  $L^\infty$-bound. The estimate of the weighted-sup norm of $\nabla u$ is obtained in a similar fashion using the pointwise gradient bounds on the Green kernel, details are omitted. In the same vein, the very last estimate (bound on the Carleson measure norm of $u$) follows from a much stronger variant which can be obtained via an integration by parts argument (testing the Poisson equation against $d(y)u$) combined with the previous $L^{\infty}$-estimate. This finishes  the proof of Lemma \ref{linearest}.
\end{proof} 

The main result of this section, pertaining to the solvability of the Dirichlet problem  \eqref{largedataeq} is the following
\begin{theorem}\label{theo42}
Let $v$ be a smooth solution of \eqref{largedataeq} subject to $v\big|_{\partial\Omega}=g\in C^1(\partial\Omega)$ and assume $v$ obeys the strict stability condition \eqref{stabcond}. Then there exists $\varepsilon>0$ such that for any $($large$)$ bounded map $f:\partial\Omega\rightarrow\mathbb{R}^m$ satisfying $\|f-g\|_{L^{\infty}(\partial\Omega)}\leq \varepsilon$, there exists a solution $u$ of the Dirichlet problem  \eqref{largedataeq}-\eqref{largedata-bc} in  $\phi+v+\W$. Moreover, this solution is unique in the ball $B^{\W}_{c_{\Omega}\varepsilon}=\{u\in \W: \|u-\phi-v\|_{\W}\leq c_{\Omega}\varepsilon\}$ for some constant $c_{\Omega}$ depending on the domain only. In particular, the solution $u$ lies in a small neighborhood of $v$, that is, $\|u-v\|_{\W}\leq \tau$ where $\tau$ depends on $\varepsilon$ and $\Omega$. Here, $\phi$ denotes the Poisson extension of $(f-g)$ to $\Omega$.  	
\end{theorem}
\begin{remark}\label{rmktheo} If $f$ is chosen large in $BMO(\partial\Omega,N)$, then the Poisson extension of $h=f-g$, $\phi_h$ is also bounded  and a similar smallness hypothesis on the $BMO$-perturbation $h$ yields the existence of a solution $u$  such that $u-v-\phi_h$ is small in $\X$. 
\end{remark}
The  stability condition in Theorem \ref{theo42}, i.e.  \eqref{stabcond} can be replaced by an invertibility condition for the linearized operator  associated to $Q_v(\cdot,\cdot)$.  Indeed, consider the bilinear form defined for any $\psi,\phi\in W^{1,2}_0(\Omega)$ with $\phi\in T_vN$ by
\begin{align*}Q_{u}(\psi,\phi):=\int_{\Omega}\bigg\{ (\nabla \psi,\nabla \phi)+\sum_{j=1}^dR^N(\phi,\partial_ju)\partial_ju\cdot\psi\bigg\} dx.\end{align*} 
Our next result shows that the conclusion of Theorem \ref{theo42} remains valid under a weaker condition.
\begin{proposition}\label{relax-invertibility} Let $\phi\in W^{1,2}_0(\Omega)$ with $\phi\in T_vN$. If condition \eqref{stabcond} is replaced by the requirement that 
\begin{equation}\label{weak-stab}
Q_{v}(\psi,\phi)=0 \hspace{0.4cm} \forall\hspace{0.01cm} \psi\in C_0^{\infty}(\Omega)\Longrightarrow \psi=0,
\end{equation}
then the conclusion of Theorem \ref{theo42} remains true. 	
\end{proposition} 

\subsection{Idea and structure of the proof of Theorem \ref{theo42}}
We discuss in this part the procedure we adopt in establishing the claims in Theorem \ref{theo42}. Once again, the plan is to perform a suitable perturbation argument in order to have a setting in which our hypotheses fit. To this end we convert the original equation into a vanishing boundary data problem. Set $h=f-g$ where $g=v\big|_{\partial\Omega}$ and denote by $\phi_{h}$ the Poisson extension of $h$ to $\Omega$. Make the ansatz $w=u-v-\phi_{h}$ and realize that $w$ solves the boundary value problem 
\begin{align*}-\Delta w&=\Gamma(v+w+\phi_{h})(\nabla (v+w+\phi_{h}),\nabla (v+w+\phi_{h}))-\Gamma(v)(\nabla v,\nabla v) \hspace{0.1cm}\text{in}\hspace{0.1cm}\Omega,\hspace{0.1cm} w\big|_{\partial\Omega}=0
\end{align*} which can be transformed into the following Dirichlet problem \begin{align}\begin{split}-\mathscr{L}_{v}w&=F(v,\phi_{h},w)\hspace{0.1cm}\text{in}\hspace{0.1cm}\Omega\\
w&=0\hspace{0.1cm}\text{on}\hspace{0.1cm}\partial\Omega\end{split}\label{mainlargeq}\end{align}
where  $-\mathscr{L}_{v}$ is the operator acting on vector-valued functions defined on $\Omega$ and given by \begin{align*}\mathscr{L}_{v}w:=\Delta w+\sum_{j=1}^{d}R^N(w,\partial_jv)\partial_jv 
\end{align*} while the nonlinearity $F(v,\phi_h,\cdot)$ reads \begin{equation*}F(v,\phi_h,w)=\Gamma(v+w+\phi_{h})(\nabla (v+w+\phi_{h}),\nabla (v+w+\phi_{h}))-\Gamma(v)(\nabla v,\nabla v)+\sum_{j=1}^{d}R^N(w,\partial_jv)\partial_jv. 
\end{equation*}  
The main focus is now on problem \eqref{mainlargeq} for it is clear that from its solvability directly flows the statement of Theorem \ref{theo42}. Assume for a moment that $-\mathscr{L}_{v}$ defined as an operator from $\W$ to $\Z$ can be inverted so that problem \eqref{mainlargeq} is reformulated as a fixed point equation $$\text{find}\hspace{0.1cm}w\in \W:\hspace{0.1cm} w=(-\mathscr{L}_{v})^{-1}\circ \widetilde{F}(v,\phi_h,w) \hspace{0.2cm}\text{in}\hspace{0.2cm}\Omega$$ where $\widetilde{F}(v,\phi_h,\cdot)$ is a suitable extension of $F(v,\phi_h,\cdot)$ which we define in subsequent lines. This allows us to avoid the geometric constraints which we treat separately. Hence, for so long as the composition $\mathcal{K}(v,\phi_h,\cdot):=(-\mathscr{L}_{v})^{-1}\circ \widetilde{F}(v,\phi_h,\cdot)$ with $v$ and $\phi_h$ as described above can be shown to be a strict contraction mapping, we are done. With other words, this amounts to saying that if altogether the following key estimate \begin{align}\label{key-est}\|w\|_{\W}\leq C\|\mathscr{L}_{v}w\|_{\Z}
\end{align} combined with the contraction property: $\exists$ $\theta_0\in(0,1)$ such that  \begin{align}\label{contratcest}\|\widetilde{F}(v,\phi_h,w_1)-\widetilde{F}(v,\phi_h,w_2)\|_{\Z}\leq \theta_0\|w_1-w_2\|_{\W}\end{align} for $w_1$ and $w_2$ in some ball of $\W$ hold true, then Theorem \ref{theo42} readily appears as a consequence of an application of the Banach fixed point theorem. Note that estimate \eqref{contratcest} only makes sense once we know that $\widetilde{F}(v,\phi_h,w)$ is an element of $\Z$ for any $w\in \W$. However, $-\mathscr{L}_{v}$ does not possess these mapping properties. In contrast, what we do know is that if the nonlinearity $\widetilde{F}(v,\phi_h,\cdot)$ maps into the Sobolev space $W^{-1,2}(\Omega)$, then the Dirichlet problem \eqref{mainlargeq} is uniquely solvable in $W_0^{1,2}(\Omega)$. Moreover, the inverse operator $(-\mathscr{L}_{v})^{-1}:W^{-1,2}(\Omega)\rightarrow W_0^{1,2}(\Omega)$ is continuous, that is, $$\|w\|_{W^{1,2}_0(\Omega)}\leq C\|\mathscr{L}_{v}w\|_{W^{-1,2}(\Omega)^m}.$$ This is a consequence of the stability condition \eqref{stabcond}, since it entails the coercivity of the continuous bilinear form associated to $-\mathscr{L}_{v}$ on $W_0^{1,2}(\Omega)$ and an application of the Lax-Milgram theorem. On the other hand, it can be checked that neither $W_0^{1,2}(\Omega)$ is a subspace of $\W$ nor $\widetilde{F}(v,\phi_h,w)$ does lie in the Sobolev space $W^{-1,2}(\Omega)$ whenever $w\in W_0^{1,2}(\Omega)$. It rather seems plausible to establish that the nonlinearity in \eqref{mainlargeq} is well-behaved with respect to the topology of the function space $\Z$. By this, we mean for every $w\in \W$ we have  \begin{align}\label{welldefbound}\widetilde{F}(v,\phi_h,w)\in \Z,\end{align}
which in turn shows that estimate \eqref{contratcest} is legitimate.
Summarizing, we shall prove that the nonlinearity $\widetilde{F}(v,\phi_h,\cdot)$ satisfies the needed properties \eqref{contratcest} and \eqref{welldefbound} -- this will constitute the first part of the proof whereas the second segment aims at showing that the operator $-\mathscr{L}_{v}$ is invertible and obeys the continuity property \eqref{key-est}. A decisive point in achieving these facts is that one has solvability of $-\mathscr{L}_{v}w=H$ in $W_0^{1,2}(\Omega)$ for $H$ in the dual space $W^{-1,2}(\Omega)$. 

To define $\widetilde{F}(v,\phi_h,\cdot)$, recall the relation between the nearest point projection map and the second fundamental form 
\begin{equation}
D^2\mathcal{P}_N(u)(V,W)=-\Gamma(u)(V,W)    
\end{equation}
whenever $u\in N$ for all $V,W\in T_uN$. Let $\mathcal{P}\in C^{\infty}(\mathbb{R}^m,\mathbb{R}^m)$ be any extension of $\mathcal{P}_N$ such that $\mathcal{P}\big|_{U_{\varrho}}=P_N$ where $U_{\varrho}:=\{x\in \mathbb{R}^m:\text{dist}(x,N)<\varrho\}$ for some $\varrho>0$ sufficiently small  and define $\widetilde{\Gamma}(z)(V,W)=D^2\mathcal{P}(z)(V,W)$ for $z\in \mathbb{R}^m$. Then 
\begin{equation*}
\widetilde{F}(v,\phi_h,w)=\widetilde{\Gamma}(v+w+\phi_{h})(\nabla (v+w+\phi_{h}),\nabla (v+w+\phi_{h}))-\widetilde{\Gamma}(v)(\nabla v,\nabla v)+\sum_{j=1}^{d}R^N(w,\partial_jv)\partial_jv.    
\end{equation*}

\begin{proof}[Proof of Theorem \ref{theo42}] 
\textbf{Part 1}. 
We state a lemma in which one quantifies the statement (\ref{welldefbound}) and also shows that (\ref{contratcest}) is indeed true.
\begin{lemma}\label{lem41}
Let $v$ falling under the scope of Theorem \ref{theo42} with $v=g$ on $\partial\Omega$ and  denote by $\phi_h$ the Poisson extension of $h=f-g$. Then the map $\widetilde{F}(v,\phi_h,\cdot)$ sends $\W$ onto $\Z$ and there exists a dimensional constant $C:=C(\Omega)$ and $K:=K(v)$ such that  \begin{align}\label{nonlinearest}\|\widetilde{F}(v,\phi_h,w)\|_{\Z}&\leq C\|w\|_{\W}(1+\|w\|_{\W})+C\|\phi_h\|_{\W}(1+\|\phi_h\|_{\W})+\\
\nonumber&\hspace{4.5cm}C(\|w\|_{\W}+\|\phi_h\|_{\W})(\|w\|^2_{\W}+\|\phi_h\|^2_{\W}+K(v)\big).\end{align}
In addition, if $\|h\|_{L^{\infty}(\Omega)}\leq \varepsilon$ for some $\varepsilon>0$ small, then there exists $\tau:=\tau(\Omega,\varepsilon)$ and $\eta\in(0,1)$ such that for all $w_1,w_2$ in the closed ball $B^{\W}_{\tau}(0)=\{w\in \W: \|w\|_{\W}\leq \tau\}$, we have $$\|\widetilde{F}(v,\phi_h,w_1)-\widetilde{F}(v,\phi_h,w_2)\|_{\Z}\leq \eta\|w_1-w_2\|_{\W}.$$
\end{lemma} 
\begin{proof}[Proof of Lemma \ref{lem41}] 
Observe that the right-hand side in \eqref{mainlargeq} can be written as		
\begin{equation*}
\widetilde{F}(v,\phi_h,w)=F_1+F_2+F_3
\end{equation*}
where
\begin{align*}
F_1&=\widetilde{\Gamma}(v)(\nabla(w+v+\phi),\nabla(w+v+\phi))-\widetilde{\Gamma}(v)(\nabla v,\nabla v), \hspace{0.2cm} F_2=\sum_{j=1}^{d}R^N(w,\partial_jv)\partial_jv,\\
F_3&=\widetilde{\Gamma}(w+v+\phi_h)(\nabla(w+v+\phi_h),\nabla(w+v+\phi_h))-
\widetilde{\Gamma}(v)(\nabla(w+\phi_h+v),\nabla(w+\phi_h+ v))
\end{align*} 
with $v$ and $\phi_h$ are as described above. Since $\widetilde{\Gamma}(\cdot,\cdot)$ and $R^N(\cdot,\cdot)\cdot$ are smooth maps, one can easily verify that
\begin{equation}\label{nonlinearity1}
|F_1|\leq c_1(|\nabla v||\nabla(w+\phi_h)|+|\nabla(w+\phi_h)|^2)
\end{equation}
where $c_1:=c_1(\|v\|_{L^{\infty}(\Omega)})$ and there exist $c_2,c_3>0$ with
\begin{equation}\label{nonlinearity2}
|F_2|\leq c_2|w||\nabla v|^2,\hspace{0.2cm}|F_3|\leq c_3|w+\phi_h||\nabla(w+v+\phi_h)|^2.
\end{equation}
Let $w\in \W$, we have \begin{align*}\|\widetilde{F}(v,\phi_h,w)\|_{\Z}&=\sup_{x\in \Omega}d(x)^2|\widetilde{F}(v,\phi_h,w)(x)|+\sup_{x\in \Omega}d(x)^{1-d}\int_{T(S)}d(y)|\widetilde{F}(v,\phi_h,w)(y)|dy\\
&=I_1+I_2.
\end{align*}
Making use of \eqref{nonlinearity1} and \eqref{nonlinearity2}, we separately estimate each of the above terms as follows.   		
\begin{align*}		
I_1	&\leq C\sup_{x\in \Omega}d(x)^2(|\nabla v||\nabla(w+\phi_h)|+|\nabla(w+\phi_h)|^2)+C\sup_{x\in \Omega}d(x)^2(|w||\nabla v|^2)+\\
&\hspace{7cm}C\sup_{x\in \Omega}d^2(x)(|w+\phi_h||\nabla (w+\phi_h+v)|^2)\\
&\leq C(v)\sup_{x\in\Omega}d(x)(|\nabla \phi_h|+|\nabla w|)+C\sup_{x\in\Omega}d^2(x)(|\nabla w|^2+|\nabla \phi_h|^2)+C(v)\|w\|_{L^{\infty}(\Omega)}+\\
&\hspace{4cm}(\|w\|_{L^{\infty}(\Omega)}+\|\phi_h\|_{L^{\infty}(\Omega)})\sup_{x\in\Omega}d^2(x)(|\nabla w|^2+|\nabla \phi_h|^2+|\nabla v|^2)\\
&\leq C\|w\|_{\W}(1+\|w\|_{\W})+C\|\phi_h\|_{\W}(1+\|\phi_h\|_{\W})+\\
&\hspace{6cm}C(\|w\|_{\W}+\|\phi_h\|_{\W})(\|w\|^2_{\W}+\|\phi_h\|^2_{\W}+C(v)\big).
\end{align*}
Taking into account the hypotheses on $v$, it follows that
\begin{align*}I_2&=\sup_{x\in \Omega}d(x)^{1-d}\int_{T(S)}d(y)|\widetilde{F}(v,\phi,w)(y)|dy\\
&\leq C\sup_{x\in \Omega}d(x)^{1-d}\int_{T(S)}d(y)(|\nabla v||\nabla(w+\phi_h)|+|\nabla(w+\phi_h)|^2)dy+\\
&\sup_{x\in \Omega}d(x)^{1-d}\int_{T(S)}d(y)|w||\nabla v|^2dy+
\sup_{x\in \Omega}d(x)^{1-d}\int_{T(S)}d(y)|w+\phi_h||\nabla(w+\phi_h+ v)|^2dy\\
&\leq C(v)\sup_{x\in \Omega}d(x)^{\frac{1-d}{2}}\big\|d(\cdot)^{1/2}|\nabla w|\big\|_{L^2(T(S))}+\sup_{x\in \Omega}d(x)^{\frac{1-d}{2}}\big\|d(\cdot)^{1/2}|\nabla \phi_h|\big\|_{L^2(T(S))}\big)+\\
&\hspace{0.5cm}\big(\sup_{x\in \Omega}d(x)^{1-d}\big\|d(\cdot)^{1/2}|\nabla w|\big\|^2_{L^2(T(S))}+\sup_{x\in \Omega}d(x)^{1-d}\big\|d(\cdot)^{1/2}|\nabla \phi_h|\big\|^2_{L^2(T(S))}\big)+\\
&\hspace{0.48cm}C(v)\|w\|_{L^{\infty}(\Omega)}+(\|w\|_{L^{\infty}(\Omega)}+\|\phi_h\|_{L^{\infty}(\Omega)})\big(\sup_{x\in \Omega}d(x)^{1-d}\big\|d(\cdot)^{1/2}|\nabla w|\big\|^2_{L^2(T(S))}+\\
&\hspace{7cm}\sup_{x\in \Omega}d(x)^{1-d}\big\|d(\cdot)^{1/2}|\nabla \phi_h|\big\|^2_{L^2(T(S))}+C(v)\big)\\
&\leq C\|w\|_{\W}(1+\|w\|_{\W})+C\|\phi_h\|_{\W}(1+\|\phi_h\|_{\W})+\\
&\hspace{6cm}C(\|w\|_{\W}+\|\phi_h\|_{\W})(\|w\|^2_{\W}+\|\phi_h\|^2_{\W}+C(v)\big).
\end{align*}Collecting the bounds on $I_1$ and $I_2$ and adding them up we conclude on the validity of \eqref{nonlinearest}. Next, we estimate $\widetilde{F}(v,\phi_h,w_1)-\widetilde{F}(v,\phi_h,w_2)$ for $w_1,w_2\in \W$ under the condition that $h=f-g$ is small in the $L^{\infty}$-norm. Write \begin{align*}\widetilde{F}(v,\phi_h,w_1)-\widetilde{F}(v,\phi_h,w_2):=A+B+C\end{align*}where 
\begin{align*}A&=\widetilde{\Gamma}(w_1+v+\phi_h)(\nabla(w_1+v+\phi_h),\nabla(w_1+v+\phi_h))-\\
&\hspace{6cm}\widetilde{\Gamma}(w_2+v+\phi_h)(\nabla(w_1+v+\phi_h),\nabla(w_1+v+\phi_h))\\
B&=\widetilde{\Gamma}(w_2+v+\phi_h)(\nabla(w_1+v+\phi_h),\nabla(w_1+v+\phi_h))-\\
&\hspace{6cm}\widetilde{\Gamma}(w_2+v+\phi_h)(\nabla(w_2+v+\phi_h),\nabla(w_2+v+\phi_h))\\
C&=\sum_{i=1}^dR^N(w_2-w_1,\partial_i v)\partial_i v\end{align*}
so that it suffices to estimate each of these quantities in $\Z$. For the same reasons as above, we have
\begin{align}\nonumber\|A\|_{\Z}&\leq C\sup_{x\in \Omega}d(x)^2|w_1-w_2||\nabla (w_1+v+\phi_h)|^2+\\
\nonumber&\hspace{5cm}C\sup_{x\in \Omega}d(x)^{1-d}\int_{T(S)}d(y)|w_1-w_2||\nabla (w_1+v+\phi_h)|^2dy\\
\nonumber&\leq C\|w_1-w_2\|_{L^{\infty}(\Omega)}\bigg(\sup_{x\in \Omega}d^2(x)|\nabla(w_1+v+\phi_h)|^2+C(v)+\\
\nonumber&\hspace{2cm}\sup_{x\in \Omega}d(x)^{1-d}\big\|d(\cdot)^{1/2}|\nabla w_1|\big\|^2_{L^2(T(S))}+\sup_{x\in \Omega}d(x)^{1-d}\big\|d(\cdot)^{1/2}|\nabla \phi_h|\big\|^2_{L^2(T(S))}\bigg)\\
&\leq C\|w_1-w_2\|_{\W}\big(\|w_1\|^2_{\W}+\|\phi_h\|^2_{\W}+C(v)\big)\label{A}
\end{align}
 Observe that $B$ can further be written as  
\begin{align*}
&\widetilde{\Gamma}(w_2+v+\phi_h)(\nabla(w_1-w_2),\nabla (w_1+v+\phi_h))+\widetilde{\Gamma}(w_2+v+\phi_h)(\nabla(w_1-w_2),\nabla (v+\phi_h))+\\
&\hspace{4cm}\widetilde{\Gamma}(w_2+v+\phi_h)(\nabla w_2,\nabla (w_1-w_2))\end{align*} from which we deduce that 
\begin{align}
\nonumber\|B\|_{\Z}&\leq C\sup_{x\in\Omega}d(x)^2(|\nabla(w_1-w_2)(x)||\nabla w_1+v+\phi_h|)+C\sup_{x\in\Omega}d(x)^2(|\nabla(w_1-w_2)||\nabla w_2|+\\
\nonumber&\hspace{1.3cm}C\sup_{x\in\Omega}d(x)^2(|\nabla (w_1-w_2(x)||\nabla(v+\phi_h)|)+\\
\nonumber&\hspace{2cm}C\sup_{x\in \Omega}d(x)^{1-d}\int_{T(S)}d(y)|\nabla(w_1-w_2)||\nabla (w_1+v+\phi_h)|dy+\\
\nonumber&\hspace{2.3cm}C\sup_{x\in \Omega}d(x)^{1-d}\int_{T(S)}d(y)|\nabla(w_1-w_2)||\nabla w_2|dy+\\
\nonumber&\hspace{3cm}C\sup_{x\in \Omega}d(x)^{1-d}\int_{T(S)}d(y)|\nabla(w_1-w_2)||\nabla (v+\phi_h)|dy\\
&\leq C\|w_1-w_2\|_{\W}(\|w_1\|_{\W}+\|w_2\|_{\W}+\|\phi_h\|_{\W}+C(v))\label{B}
\end{align}
where we have applied H\"{o}lder's inequality to estimate the integral terms. Finally, we have 
\begin{align}\nonumber\|C\|_{\Z}&\leq C\sup_{x\in\Omega}d(x)^2|\nabla v|^2|w_1-w_2|+C\sup_{x\in \Omega}d(x)^{1-d}\int_{T(S)}d(y)|\nabla v|^2|w_1-w_2|dy\\
&\leq CC(v)\|w_1-w_2\|_{\W}.\label{C}
\end{align}
The generic constant appearing in \eqref{A}, \eqref{B} and \eqref{C} depends on $\mathrm{diam}$ $\Omega$.
Recall that $\phi_h$ is small in $\W$ since $h$ is small in $L^{\infty}(\partial\Omega)$ from Lemma \ref{linearest}. For $w_1$ and $w_2$ in the closed ball $B^{\W}_{\tau}(0)$ of $\W$, we deduce in view of the above bounds on $A$, $B$ and $C$ that the second part of Lemma \ref{lem41} holds true. The proof of Lemma \ref{lem41} is now complete.
\end{proof}
\textbf{Part 2}. Here we prove the invertibility of the operator $-\mathscr{L}_{v}$ together with \eqref{key-est}. Let us introduce the operator $L:=\Delta +\ell$ acting on $\mathbb{R}^m$-valued functions defined on $\Omega$ where $\ell$ is a smooth linear map. Consider the zero data Dirichlet boundary value problem \begin{align}\begin{split}-Lw(x)&=F(x),\hspace{0.2cm}x\in \Omega\\ w(x)&=0,\hspace{0.2cm}x\in \partial\Omega\end{split}\label{iterationeq}
\end{align}
where $F\in\mathbf{Z}$. Remark that the operator $\mathscr{L}_{v}$ has the form of $L$ with \[\ell(w)=\sum_{j=1}^{d}R^N(w,\partial_jv)\partial_jv.\] 
Now if $w^0$ is a solution of the Poisson equation $-\Delta w^0=F$ in $\Omega$ with zero data at the boundary, then $w^1=w-w^0$ solves the Dirichlet problem \begin{align}\begin{split}
-Lw^1&=F_0:=\sum_{j=1}^{d}R^N(w_0,\partial_jv)\partial_jv\hspace{0.2cm}\text{in}\hspace{0.2cm}\Omega\\
w^1\big|_{\partial\Omega}&=0.\end{split}\label{w1eq}
\end{align}
At this point, we only need to show that $w^1$ belongs to $\W$ with a corresponding \textquotedblleft good'' estimate since the solution $w^0$ is well understood by now due to Lemma \ref{linearest}. A first step towards this is the following\begin{claim}$F_0$ belongs to $\Z\cap W^{-1,2}(\Omega)$.\end{claim}  
Let us momentarily defer the proof of this claim and observe that it implies the existence of a unique $w^1\in W_0^{1,2}(\Omega)$ solving \eqref{w1eq}. The extra information $F_0\in \Z$ will enable us to improve the regularity of $w^1$ via an iterative scheme. Set $F_1=\sum_{j=1}^{d}R^N(w^1,\partial_jv)\partial_jv$ and let $w^2$ be such that $-\Delta w^2=F_0$ in $\Omega$. One can easily verify that $F_1\in L^2(\Omega)$ and that $w^3=w^1-w^2$ is a solution to the problem $-\Delta w^3=F_1$ with zero data on $\partial\Omega$. This implies by elliptic regularity theory $w^3\in W^{2,2}(\Omega)$. Iterating this procedure, we eventually find $w^1\in\W$. Hence, $w=w^1+w^0\in \W$ and $\|w\|_{\W}\leq C\|F\|_{\Z}$.  
This, in concert with part 1, proves \eqref{key-est}. By the method of continuity,  $-\mathscr{L}_{v}$ is invertible. Moreover, for $w_1,w_2\in B_{\tau}^{\mathbf{W}}(0)$, we have in view of \eqref{key-est} and using Lemma \ref{lem41} from part 1, \begin{align*}\|\mathcal{K}(v,\phi,w_1)-\mathcal{K}(v,\phi,w_2)\|_{\mathbf{W}}&\leq C\|\widetilde{F}(v,\phi,w_1-w_2)\|_{\Z}\\
	&\leq \varepsilon\|w_1-w_2\|_{\W} 
\end{align*}where $\varepsilon=\varepsilon(\tau)$ can be made small if $\tau$ is sufficiently small. Thus, $\mathcal{K}(v,\phi,\cdot)$ is a strict contraction mapping. We now verify that $F_0\in \Z\cap W^{-1,2}(\Omega)$.
\begin{align*}\|F_0\|_{\Z}&=\sup_{x\in \Omega}d(x)^2|F_0(x)|+\sup_{x\in \Omega}d(x)^{1-d}\int_{T(S)}d(y)|F_0(y)|dy\\
&=I+II.
\end{align*}
Using the smoothness of $R^N$ and the fact that $w^0\in\W $, one finds that
\begin{align*}
I&\leq \sup_{x\in \Omega}d(x)^2\bigg|\sum_{j=1}^{d}R^N(w^0,\partial_jv)\partial_jv\bigg|\leq C(v)\|w^0\|_{\W}\\
&\leq C(v)\|F\|_{\Z}.
\end{align*}
For the solid integral, we have
\begin{align*}II&\leq \sup_{x\in \Omega}d(x)^{1-d}\int_{T(S)}d(y)\bigg|\sum_{j=1}^{d}R^N(w^0,\partial_jv)\partial_jv\bigg|dy\\
&\leq C(v)\|w^0\|_{L^{\infty}(\Omega)}\\
&\leq C(v)\|w^0\|_{\W}\leq C(v)\|F\|_{\Z}.
\end{align*}
Now we establish that $F_0\in W^{-1,2}(\Omega)$. Let $\varphi\in W^{1,2}_0(\Omega)$, using the integration by parts formula and H\"{o}lder's inequality, it follows that \begin{align*}
\bigg|\left\langle F_0,\varphi \right\rangle_{W^{-1,2},W_0^{1,2}}\bigg|&=\bigg|\int_{\Omega}(\sum_{j=1}^{d}R^N(w^0,\partial_jv)\partial_jv)\cdot\varphi dx\bigg|\\
&\leq C(v)\|w^0\|_{\W}\|\varphi\|_{L^2(\Omega)}\leq C(v)\|F\|_{\Z}\|\varphi\|_{W^{1,2}(\Omega)}.
\end{align*} 
To conclude the proof of Theorem \ref{theo42}, one needs to show that $w+\phi_h+v\in N$. Since $v\in N$, it holds that
\begin{align*}
\text{dist}(w+\phi_h+v,N)&\leq C\|w+\phi_h\|_{L^{\infty}(\Omega)}\\
&\leq C(\|w\|_{\W}+\|\phi_h\|_{\W})\leq C\varepsilon
\end{align*}
for $\varepsilon>0$ small. Thus $w+\phi_h+v\in U_{\varepsilon'}$, $\varepsilon'=C\varepsilon$ (with $U_{\varepsilon'}$ as defined above) so that one can define $\Upsilon_N(w+\phi_h+v)=w+\phi_h+v-\mathcal{P}_N(w+\phi_h+v)$. It is clear that $\Upsilon_N(w+\phi_h+v)\big|_{\partial\Omega}=0$ because $(w+\phi_h+v)\big|_{\partial\Omega}\in N$. Moreover, similar calculations to those performed in Section \ref{S3} reveals that $\Upsilon_N(w+\phi_h+v)$ is subharmonic in $\Omega$. As a consequence, $\Upsilon_N(w+\phi_h+v)=0$ in $\Omega$. This achieves the Proof of Theorem \ref{theo42}.
\end{proof} 
\begin{proof}[Proof of Proposition \ref{relax-invertibility}]
We first solve the inhomogeneous linear problem $-\mathscr{L}_vw=H$ in $W^{1,2}_0(\Omega)$ for $H$ in $\mathbf{Z}$. The remaining bit of the proof will just be a reprise of the argument in part 2 of the proof of Theorem \ref{theo42}. Again, write $\mathscr{L}_v:=\Delta+\ell$ with  $\ell(\cdot)=\sum_{j=1}^dR^N(\cdot,\partial_jv)\partial_jv$ and set $K=(-\Delta)^{-1}\circ \ell$. The operator $\ell$ is bounded and compact from $W^{1,2}_0(\Omega)$ to $L^2(\Omega)$. On the other hand, the inverse Laplacian $(-\Delta)^{-1}:L^2(\Omega)\rightarrow L^2(\Omega)$ is also bounded and compact so that $K$ can be realized as a linear bounded compact operator from $W^{1,2}_0(\Omega)$ to $L^2(\Omega)$ and $\widetilde{H}=(-\Delta)^{-1}H\in L^2(\Omega)$ since $(-\Delta)^{-1}$ maps continuously $\mathbf{Z}$ into $\mathbf{W}$. Our problem therefore reduces to that of solving 
	\[
	\begin{cases}
	w+Kw=\widetilde{H} \hspace{0.1cm}\text{in}\hspace{0.1cm} \Omega\\
	w\rvert_{\partial\Omega}=0.
	\end{cases}
	\]
By virtue of the hypothesis in Proposition \ref{relax-invertibility}, the trivial solution is  the only solution of $w+Kw=0$ in $\Omega$ with vanishing boundary data. Hence, existence of a unique solution for the above problem in $W^{1,2}_0(\Omega)$ is a consequence of the Fredholm alternative. The conclusion then follows from the proof of Theorem \ref{theo42} (see part 2) as previously mentioned.     	
\end{proof}

\bibliographystyle{acm}

\end{document}